\documentclass[10pt,twocolumn]{IEEEtran}

\usepackage[noadjust]{cite}
\usepackage{subfigure}
\usepackage{epsfig}
\usepackage{courier}
\usepackage{amssymb,amsmath,color,graphicx, amsbsy}
\usepackage{algorithm}
\usepackage{algorithmic}
\usepackage{amsfonts}
\usepackage{epstopdf}
\usepackage{booktabs, multicol, multirow}
\usepackage{pifont}

\usepackage{supertabular,booktabs}

\usepackage{caption}  

\usepackage{soul}
\usepackage{subfig} 
\usepackage{array}  

\usepackage[dvipsnames]{xcolor}  

\DeclareCaptionFont{MyFigureFont}{\footnotesize}

\newtheorem{theorem}{Theorem}

\newtheorem{lemma}{Lemma}

\makeatletter
\def\mod@estimate@lineht{%
  \ST@lineht=\arraystretch \baslineskp
  \global\advance\ST@lineht by 1\p@
  \ST@stretchht\ST@lineht\advance\ST@stretchht-\baslineskp
  \ifdim\ST@stretchht<\z@\ST@stretchht\z@\fi
  \ST@trace\tw@{Average line height: \the\ST@lineht}%
  \ST@trace\tw@{Stretched line height: \the\ST@stretchht}%
}
\newenvironment{strictsupertabular}
  {\let\estimate@lineht\mod@estimate@lineht\supertabular}
  {\endsupertabular}
\makeatother

\usepackage{etoolbox} 

\let\mybibitem\bibitem
\renewcommand{\bibitem}[1]{%
  \ifstrequal{#1}{nature}
    {\color{blue}\mybibitem{#1}}
    {\color{black}\mybibitem{#1}}%
}

\begin{document}

\title{Strategic Bidding for  Producers  in Nodal Electricity Markets: A Convex Relaxation Approach}

\author{Mahdi Ghamkhari,~\IEEEmembership{Student Member,~IEEE}, Ashkan Sadeghi-Mobarakeh,~\IEEEmembership{Student Member,~IEEE}   \\   and Hamed Mohsenian-Rad,~\IEEEmembership{Senior Member,~IEEE} \thanks{ }
\vspace{-0.35cm}
 \thanks{The authors are with the Department of Electrical Engineering, University of California, Riverside, CA, USA, e-mails: \{ghamkhari, asade004, hamed\}@ee.ucr.edu. This work was supported by NSF grants 1253516, 1307756, and 1319798. The corresponding author is H. Mohsenian-Rad.}}

\maketitle

\begin{abstract}

Strategic bidding problems in electricity markets are widely studied in power systems, often by formulating complex bi-level optimization problems that are hard to solve. The state-of-the-art approach to solve such problems is to reformulate them as mixed-integer linear programs (MILPs). However, the computational time of such MILP reformulations grows dramatically, once the network size increases, scheduling horizon increases, or randomness is taken into consideration. In this paper, we take a fundamentally different approach and propose effective and customized convex programming tools to solve the strategic bidding problem for producers in nodal electricity markets. Our approach is inspired by the Schmudgen's  Positivstellensatz Theorem in semi-algebraic geometry; but then we go through several steps based upon both convex optimization and mixed-integer programming   that results in obtaining close to optimal    bidding solutions, as evidenced by several numerical case studies, besides having a huge advantage on reducing computation time. While the computation time of the state-of-the-art MILP approach grows exponentially when we increase the scheduling horizon or the number of random scenarios, the computation time of our approach increases rather linearly.


\vspace{0.25cm}

\textbf{\emph{Keywords}}: Nodal electricity market,  strategic bidding, equilibrium constraints, convex optimization, computation time.

\end{abstract}


\vspace{-.1cm}
\section*{Nomenclature}

\noindent

\begin{strictsupertabular}{ll}

$\mathbb{R}$, $\mathbb{R}^+$&Set of real and non-negative real  numbers \\


$\mathbb{S}$ & Set of symmetric matrices \\

$\mathcal{N}$ & Set of nodes in power grid in arbitrary order \\

$\mathcal{D}$ & Set of demand nodes in ascending order \\

$\mathcal{G}$ & Set of generation nodes in ascending order \\

$\mathcal{S}$ & Subset of strategic generation nodes in set $\mathcal{G}$ \\

$\mathcal{L}$ & Set of transmission lines, in  arbitrary order \\

$k$ & Index for random scenarios \\

$[t]$ & Hourly time slots\\

$T$ & Number of hourly time slots\\

$K$ &Number of random scenarios \\

$P_G$ & Vector of power generations \\

$P_D$ &  Vector of demands  \\

$\theta$ & Vector of phase angels of power grid\\

$\lambda$ & Vector  of locational marginal prices \\

$\sigma$, $\delta$, $\zeta$, & Vectors  of dual variables corresponding  \\

$\xi$, $\phi$, $\psi$ & to  inequalities in economic dispatch problem \\

$A$& Bus-line incidence matrix\\

$B_G$& Generator-bus incidence matrix\\

$B_D$ &   Demand-bus incidence matrix\\

$B_S$&  Strategic generators to  generators  incidence matrix  \\

$V$ &  Diagonal matrix of transmission lines  reactance  \\

$a$ & Vector of  energy  price bid of generators \\

$b$& Vector of demand price bid  of loads \\

$c$ &Vector of cost parameter of  strategic generators   \\
$C$ &  Vector of line capacities \\

$P_G^{\text{min}}, P_G^{\text{max}}$ & Vector of minimum and maximum generation\\

$P_D^{\text{min}}, P_D^{\text{max}}$ & Vector of minimum and maximum demand\\

$\Gamma$ & Ramp constraint parameter \\

$\mathbf{0}$ & A  column vector or a matrix with zero entries \\

$x$ & Column vector of all variables in (\ref{formula:MPEC_KKT}) \\

$n$ & Length of the vector $x$ \\

$F,Q$ & Symmetric matrices of parameters  in $\mathbb{S}^{n}$ \\
$f,p,v, q, d$ & Vectors of parameters  in $\mathbb{R}^n$\\


$r, O, \bar{x}$ & Defined in (\ref{formula:Matlab_Command}) \\

 $*$   &  Point-wise production of two vectors  \\

$\text{Rank}(\cdot)$ & Rank of a matrix \\

${(\cdot)^T}$ & Transpose of a vector or a matrix  \\




$tr(\cdot)$& Trace of a matrix\\

$\succeq$ &  Matrix inequality\\
$ i,j$  &  Indices for   $I$ linear inequalities in (\ref{formula:QPQC}),  $i,j\leq I$\\
$z$  &  Index for   $Z$ quadratic equalities in (\ref{formula:QPQC}), $z\leq Z$\\

$m$  &   Index for   $M$ linear equalities in (\ref{formula:QPQC}), $m\leq M$\\

$l$ & Index for $n$ elements of a vector in $\mathbb{R}^n$, $l\leq n$\\
$e_l$ & $l$th element of the standard basis for $\mathbb{R}^n$ space  \\

\end{strictsupertabular}

\vspace{.1cm}
\section{Introduction}

Strategic bidding plays a central role in wholesale electricity markets, where  market participants seek to choose their  bids to the day-ahead and/or real-time markets so as to maximize their  profits. Strategic bidding  in electricity markets has been extensively  studied previously, e.g.,  for producers \cite{ruiz2009pool,MILP_Example_1,MILP_Example_2, Market_Power_2}, consumers \cite{Lagrangian_Relaxation_2,DR1,DR2}, and energy storage units \cite{ESS1,Mohsenian_MPEC_2015,ESS3}.

The literature on strategic bidding is often categorized based on whether the market participant is small and price-taker \cite{DR1,ESS1, Ladurantaye}, or large and price-maker \cite{ESS3,Market_Power_2,DR2,ruiz2009pool,MILP_Example_1,MILP_Example_2,Mohsenian_MPEC_2015,Lagrangian_Relaxation_2}. The focus in this paper is on the latter, where the details on how the  market operates are \emph{explicitly} considered in formulating the strategic bidding problem. Accordingly, the strategic bidding problem  is formulated as a \emph{bi-level} program, where the lower  level problem constitutes  the  economic dispatch problem that is   solved by the independent system operator (ISO) in order  to minimize the cost of electricity dispatch and to set the  market prices. Following the common approach in the electricity market literature,  the strategic bidding problem is then reformulated  as  a single \emph{mathematical program with equilibrium constraints} (MPEC),  see \cite{ruiz2009pool},   \cite{Hobbs2000,A7,A8}.

A wholesale market   offering strategy is proposed in \cite{A7} for a wind power producer
with market power, which participates in the day-ahead market as a price-maker, and in the balancing market as a deviator. Optimal bidding for a large consumer is formulated in \cite{A3} as an MPEC problem. MPEC formulation is also used in \cite{A8} for optimal strategic bidding of a regulation resource in the performance-based regulation market considering the system dynamics.
 The 
preventive maintenance scheduling of power transmission lines within a yearly time
framework using a bi-level optimization approach was studied in \cite{A9}.  In  \cite{A10},
a vulnerability analysis of an electric grid under disruptive threat is formulated as a 
bi-level optimization.  Finally, in \cite{Hobbs2000},  strategic gaming in electricity markets  was   analyzed  using 
 an MPEC formulation.

 The MPEC problem formulations that appear in power systems are often difficult to solve. The difficulty arises due to the necessary use of bilinear terms that create non-convex objective function and constraints. The common approach to solve such problems  is to transform them into mixed integer linear programs  (MILPs), e.g., see 
 \cite{ruiz2009pool,MILP_Example_1,MILP_Example_2,Mohsenian_MPEC_2015,DR2,ESS3}.

While the MILP reformulations of strategic-bidding problems are popular in the power systems community, such reformulations are prone to \emph{major computational challenges}. Specifically, the computational time often increases dramatically, once the network size grows, scheduling horizon increases, or randomness is taken into consideration. For example, for one of our case studies with  10  random scenarios, the MILP approach in \cite{ruiz2009pool} did not converge even after letting it run for about three days, see Section \ref{sec:SingleTimeSlot} for details.   

To tackle the aformentioned computational  challenges,  some attempts \emph{with little success} have  been made recently to solve the strategic bidding problems in power systems using convex optimization techniques. In particular, in \cite{Fampa} and \cite{haghighat2014strategic}, the authors used \emph{semidefinite relaxation} and \emph{lift-and-project linear relaxation} to solve the MPEC problems in electricity markets. However, in both cases, the performance was often poor with respect to not only optimality but also computation time. Moreover, no clear recovery method was proposed to guarantee obtaining a feasible solution of the original MPEC problem. Finally, only small MPEC problems were discussed.
%
%

Therefore, to the best of our knowledge, it is fair to say that solving the strategic bidding problems in wholesale electricity markets using convex programming is still an open problem and no reliable  and scalable  solution approach currently exists to address the  relatively large and hence practically relevant problems. Accordingly, our goal in this paper is to tackle this open problem.  Without loss of generality, we focus on the case of  strategic bidding for \emph{producers}. The main technical contributions in this paper can be summarized as follows:
\begin{itemize}

\item We take a fundamentally different approach from \cite{ruiz2009pool,MILP_Example_1,MILP_Example_2} and \cite{Fampa, haghighat2014strategic}, and propose innovative and effective convex programming tools to solve the strategic bidding problem for producers in nodal electricity markets, where our approach is \emph{customized} to exploit the main characteristics of such problems. Our proposed solution method is accurate, reliable, and computationally tractable in solving the strategic bidding problems in power systems.

\vspace{0.1cm}

\item Our approach is initially inspired by the  Schmudgen's  Positivstellensatz Theorem \cite[Theorem 3.16]{Laurent_Book_2009} \cite[Section 4.3]{SOS_Manual} in semi-algebraic geometry; but then we go through several steps based upon both convex optimization and mixed-integer programming in order to develop an algorithm, Algorithm 1, that is guaranteed to give a feasible and very close-to-optimal solution to the original MPEC problem, besides having a huge advantage on reducing computation time.

\vspace{0.1cm}


\item We compare the optimality and the computation time of our proposed approach and that of the MILP approach in \cite{ruiz2009pool} for the case of a market over the IEEE 30 bus test system. While the computation time of the MILP approach in \cite{ruiz2009pool} increases \emph{exponentially} when we increase the scheduling horizon or the number of random scenarios, the computation time of our proposed approach increases rather \emph{linearly}. Interestingly, the average optimality of the solution from our proposed approach is $99\%$ or higher.

\end{itemize}

It is worth pointing out that the state-of-the-art  polynomial optimization problem relaxations that are formulated based on   Schmudgen's Positivestellensatz \cite[Theorem 3.16]{Laurent_Book_2009} and  Lasserre's sum-of-squares \cite{Jeyakumar2014,Nie2012}   tend to provide  tight upper bounds for  the intended non-convex optimization problems only when we significantly increase    the order of  added coefficients or  polynomials. Accordingly, in both cases, we often face convex but very large optimization problems for any descent size problem, which makes the resulting convex relaxation approach of little interest in practice.  In contrast,  in this paper, we use  Schmudgen's Positivestellensatz but not  Lasserre's sum-of-squares  method,  because we are able to build upon  it a new methodology, combined with a heuristic algorithm, which results in obtaining very  close to optimal bidding solutions within a reasonable computational time.

The proposed approach in this paper can be applied to the other MPEC problems in electricity markets, e.g., to find optimal bids for large energy storage units \cite{Mohsenian_MPEC_2015}, or to tackle the strategic generation investment problem for producers \cite{MILP_Example_1}.

\vspace{.1cm}
\section{Problem Statement}

Consider a strategic price-maker generation firm that bids in a day-ahead nodal electricity market. Once the bids from all market participants  are collected, the ISO solves an economic dispatch problem, which is presented below in vector-format, in order  to determine  the clearing market price and   the  energy reward to each producer \cite[Appendix C]{Conejo_Book}, \cite{Y_Fu_2006}:
\begin{alignat}{8}
 &\underset{P_G,P_D,\theta}{\textbf{minimize}} \ \ {a}^TP_G-{b}^TP_D \label{formula:Dispatch}\\
&\textbf{subject to} \nonumber\\
&B_GP_G-B_DP_D-AV^{-1}A^T\theta =0 :\lambda \   \label{formula:KKL}\\
&P_G -P_G^{\text{min}}  \geq \mathbf{0}  : \sigma \label{formula:P_G_Limit_Down}\\
&P_G^{\text{max}} - P_G \geq \mathbf{0} :\delta  \label{formula:P_G_Limit_UP}\\
& P_D -P_D^{\text{min}} \geq \mathbf{0}: \zeta   \label{formula:P_D_Limit_Down}\\
&P_D^{\text{max}}- P_D \geq \mathbf{0} : \xi\label{formula:P_D_Limit_UP}\\
&V^{-1}A^T\theta + C  \geq \mathbf{0}  :\phi \label{formula:Capacity_Down} \\
&C - V^{-1}A^T\theta \geq \mathbf{0} :\psi,  \label{formula:Capacity_up}
\end{alignat}
where the notations are explained in the Nomenclature.  The  vector of power flows on all  transmission lines  is modeled here as   ${V}^{-1}A^T\theta$. The variable after each colon  in (\ref{formula:KKL})-(\ref{formula:Capacity_up}) shows the dual variable  corresponding to each constraint. Constraint (\ref{formula:KKL}) enforces the power balance  and its dual variable is  the market clearing price. Constraints (\ref{formula:P_G_Limit_Down})-(\ref{formula:P_D_Limit_UP}) enforce the generations and loads to operate  within their  limits. Furthermore, constraints (\ref{formula:Capacity_Down})-(\ref{formula:Capacity_up}) enforce the  capacity  for transmission lines.

Note that, for the ease of discussions, the problem formulation in (\ref{formula:Dispatch})-(\ref{formula:Capacity_up}) is for economic bidding over a single hour. The case for multiple hours  is discussed later in Section \ref{sec:Multi}.

\vspace{.1cm}
\subsection{Bi-level Problem Formulation}\label{sec:bilevel}
The bidding  problem  for the strategic generation firm of interest  can be formulated  as a \emph{bi-level} program  \cite{ruiz2009pool, MILP_Example_1}:
\begin{equation}\label{formula:MPEC}
\begin{aligned}
 & \!\! \underset{\begin{subarray}{c}B_S a,P_G,P_D,\phi\\ \theta, \lambda,\sigma,\delta,\zeta,\xi,\psi\end{subarray}}{\textbf{maximize}}  && \!\!  {\lambda}^TB_G {B_S}^T B_S  P_G-  c^TB_S P_G  &&&& \\
&\textbf{subject to}   && \!\!\!\left( \begin{subarray}{c}P_G,P_D\\ \phi, \theta, \lambda, \sigma \\  \delta, \zeta,\xi,\psi\end{subarray}\right) =  &&  \!\!\!\!\!\!\!\!\!\!\!\!\!\!\!\!\!\!\!\!\!\!\!\!\!\!\!\!\!\!\!\!\!\!\!\!\!\!\!\!\!\!\!\!\! \textbf{argmin} \; && \!\!\!\!\!\!\!\!\!\!\!\!\!\!\!\!\!\!\!\!\! a^TP_G-b^TP_D\\
&&&&&\!\!\!\!\!\!\!\!\!\!\!\!\!\!\!\!\!\!\!\!\!\!\!\!\!\!\!\!\!\!\!\!\!\!\!\!\!\!\!\!\!\!\!\!\! \textbf{subject to}  \; &&\!\!\!\!\!\!\!\!\!\!\!\!\!\!\!\!\!\!\!\!\!(\ref{formula:KKL})-(\ref{formula:Capacity_up}).
 \end{aligned}
\end{equation}
  The two terms in the objective function in (\ref{formula:MPEC}) denote the \emph{total generation revenue} and \emph{the total generation cost} for the strategic generation firm of interest,  respectively. The upper-level problem in (\ref{formula:MPEC}) constitutes the profit maximization problem that the strategic generation firm seeks to solve. The lower-level problem in (\ref{formula:MPEC}) constitutes the   economic dispatch problem that the ISO must solve, before  the profit  of the generation firm can be calculated at the upper-level problem.  Note that, the  optimization variables in problem (\ref{formula:MPEC}) include  the vector of price bids for strategic generators, which is represented here  as $B_S a$,  where $B_S\in\mathbb{R}^{|\mathcal{S}|\times|\mathcal{G}|} $ is the incidence matrix for the vector of  strategic generators to the vector of all generators, and $a$  is the vector of price bids for all generators.  The elements of  vector $a$ that  belong to the set of   non-strategic   generators are taken as  parameters in problem (\ref{formula:MPEC}). 
  
 In formulating problem (\ref{formula:MPEC}), we followed the same assumption as in \cite[Section 4]{Garcia2016} and \cite{Dempe2013}   in the sense that if there exist multiple solutions for problem (\ref{formula:Dispatch})-(\ref{formula:Capacity_up}), then the solution   that is most profitable to  the firm is considered.     

\vspace{.1cm}
\subsection{MPEC Problem Reformulation}\label{sec:MPECfor}
The lower-level problem in (\ref{formula:MPEC}) is a linear program. Therefore, its corresponding   Karush-Kuhn-Tucker (KKT) optimality conditions  are  both necessary and  sufficient  \cite[Section 5.5.3]{Boyd}.  They comprise (\ref{formula:KKL})-(\ref{formula:Capacity_up}),  and  the following constraints:
\begin{align}
&\begin{aligned}
&a- B_G\:  \lambda-\sigma+\delta=0
\end{aligned}& \label{formula:Equal_G} \\
&b- B_G\:   \lambda+\zeta-\xi=0  &  \label{formula:Equal_D}\\
&\begin{aligned}
&AV^{-1}(A^T\lambda+ \psi-\phi)=0
\end{aligned} &  \label{formula:LMP}\\
&\sigma * (P_G-P_G^{\text{min}})=0   &  \label{formula:Sigma} \\
&\delta * (P_G^{\text{max}}-P_G)=0  &  \\
&\zeta * (P_D-P_D^{\text{min}})=0   &    \\
&\xi * (P_D^{\text{max}}-P_D)=0   &   \\
&\phi * (C+V^{-1}A^T\theta)=0   &   \\
&\psi * (C-V^{-1}A^T\theta)=0   &  \label{formula:psi}
\end{align}
\begin{align}
&\begin{aligned}
\delta\geq 0,\
\xi\geq 0, \
\psi\geq 0
\end{aligned} &  \label{formula:Down} \\
&\begin{aligned}
\sigma\geq 0,\
\zeta\geq 0,\
\phi\geq 0.
\end{aligned} &  \label{formula:Up}
\end{align}
Once we replace the  lower-level problem     with its equivalent KKT conditions, the  bi-level strategic bidding problem in (\ref{formula:MPEC}) takes the form of a   standard MPEC problem as follows:
\begin{equation}\label{formula:MPEC_KKT}
\begin{aligned}
& \!\! \underset{\begin{subarray}{c}B_S a,P_G,P_D,\phi\\ \theta, \lambda,\sigma,\delta,\zeta,\xi,\psi\end{subarray}}{\textbf{maximize}} \ {\lambda}^TB_G {B_S}^T B_S  P_G-  c^TB_S P_G  \\
&\textbf{subject to}\   (\ref{formula:KKL})-(\ref{formula:Capacity_up}) \ \text{and} \ (\ref{formula:Equal_G})-(\ref{formula:Up}).
\end{aligned}
\end{equation}
 Problem (\ref{formula:MPEC_KKT}) is non-convex and hard to solve.  Non-convexity is due to the bilinear terms, both in the complimentary slackness constraints   (\ref{formula:Sigma})-(\ref{formula:psi}) and in the first term in the  objective function.  
For the rest of this paper, we seek to solve problem (\ref{formula:MPEC_KKT}) in an \emph{accurate} yet  \emph{computationally tractable} fashion. %

  We assume that the economic dispatch problem in (\ref{formula:Dispatch})-(\ref{formula:Capacity_up}) is  always  feasible \cite{zhao2015distributed,li2014connecting}. Note that, since   this problem is a linear program, it always satisfies the  slater's constraints qualifications conditions \cite[Section 5.2.3]{Boyd}. Therefore, there  always  exists a solution for the  KKT conditions of the economic dispatch problem. That is, the set of constraints in (\ref{formula:KKL})-(\ref{formula:Capacity_up}) and (\ref{formula:Equal_G})-(\ref{formula:Up})  is always feasible. Therefore, the  MPEC problem in (\ref{formula:MPEC_KKT}) always has a  feasible solution.

\vspace{.1cm}
\section{Solution Method}\label{sec:Solution}

\vspace{.1cm}
\subsection{Fundamental Convex Relaxation Approach}\label{sec:SolutionA}
 The common approach to solve  problem (\ref{formula:MPEC_KKT}) is to reformulate it as a mixed  integer linear program, e.g., see  \cite{ruiz2009pool,MILP_Example_1,MILP_Example_2,Mohsenian_MPEC_2015}.  However, the computation time of solving such MILP reformulation grows exponentially  as the size of problem (\ref{formula:MPEC_KKT}) increases \cite{ruiz2009pool}. Therefore, in  this section, we present an alternative approach to solve problem (\ref{formula:MPEC_KKT}) based on  convex  optimization,  where computation time grows  linearly.   We concatenated all the optimization variables in problem (\ref{formula:MPEC_KKT}) into a single optimization vector as follows:
 \begin{equation}
\begin{aligned}
x\triangleq{[{(B_S a)}^T  P_G^T\;\! P_D^T\;\!  \lambda^T\;\! \sigma^T \;\!\delta^T\;\! \zeta^T \;\!\xi^T \;\! \phi^T\;\! \psi^T\;\! \theta^T]}^T.
\end{aligned}
\end{equation}
Let $n$ denotes the length of vector $x$.
 First, we represent   problem (\ref{formula:MPEC_KKT}) in its vector form   as follows \cite[Section 4.4]{Boyd}:
\begin{equation}\label{formula:QPQC}
\begin{aligned}
&\underset{x}{\textbf{maximize} }&    &x^TFx+2f^Tx &                          &                       \\
&\textbf{subject to}                      &    &p_i^Tx+p_{i0}\geq 0&                  & \forall i      \\
&                                                  &    &v_m^Tx+v_{m0}=0&                    &\forall m  \\
&                                                  &    & x^TQ_zx+2q_z^Tx=0&   &\forall z   ,
\end{aligned}
\end{equation}
  where $x\in\mathbb{R}^n$ is the column vector of all decision variables in problem (\ref{formula:MPEC_KKT}).  Here, $F$ and $f$ are derived from the objective function in (\ref{formula:MPEC}); $p_i$ and $p_{i0}, \forall i$, are derived from the linear inequality constraints in (\ref{formula:P_G_Limit_Down})-(\ref{formula:Capacity_up}), (\ref{formula:Down}), (\ref{formula:Up}); $v_m$ and $v_{m0}, \forall m$, are derived from the linear equality constraints in (\ref{formula:KKL}), (\ref{formula:Equal_G})-(\ref{formula:LMP}); and $Q_z$ and $q_z, \forall z$, are derived from the quadratic equality constraints in (\ref{formula:Sigma})-(\ref{formula:psi}). Since all quadratic equality constraints  are due to \emph{complimentary  slackness}, we can write
\begin{equation}\label{formula:Slack_Feature}
Q_z=d_z q_z^T \ \ \ \forall z,
\end{equation}
 where  $d_z$, $\forall z$, is derived from (\ref{formula:Sigma})-(\ref{formula:psi}). We will use (\ref{formula:Slack_Feature}) later in Section \ref{sec:Recovery}.   Problem  (\ref{formula:QPQC}) is always feasible, since it is a  reformulation of problem (\ref{formula:MPEC_KKT}), see   Section \ref{sec:MPECfor}.

 Problem  (\ref{formula:QPQC})  is a quadratically-constrained quadratic program (QCQP).   Following the  analysis in \cite[Theorem 3.16]{Laurent_Book_2009},  we propose the following \emph{relaxation} of problem (\ref{formula:QPQC}):
 \begin{equation}\label{formula:QPQC_Relaxed}
\begin{aligned}
&\underset{\begin{subarray}{c}\Lambda,\alpha_i,\varrho_{ij}\\ \beta_z,h_m,h_{m0} \end{subarray}}{\textbf{minimize} }       &&  \Lambda  \\
&\textbf{subject to} \!\!\! &&   \Lambda -x^TFx-2f^Tx -\sum_{i=1}^I\alpha_i(p_i^Tx+p_{i0})- \\
&&&\sum_{i=1}^{I}\sum_{j=1}^I \varrho_{ij} (p_i^Tx+p_{i0})(p_j^Tx+p_{j0})  - \\
&&& \sum_{m=1}^M( h_m^T x+h_{m0})(v_m^Tx+v_{m0}) -\\
&&&\sum_{z=1}^Z \beta_z(x^TQ_zx+2q_z^Tx)  \geq 0 \ \ \ \forall x\in\mathbb{R}^{n},
\end{aligned}
\end{equation}
 where $\Lambda\in\mathbb{R}$, $\alpha_i\in\mathbb{R}^+$ and $\varrho_{ij}\in\mathbb{R}^+$ $\forall i$ and $\forall j$, $ \beta_z\in\mathbb{R} $  $  \forall z $, $h_m\in\mathbb{R}^n$ $\forall m$, and $h_{m0}\in\mathbb{R}$ $\forall m$.   We shall point out four key properties of problem (\ref{formula:QPQC_Relaxed}).  First,   $x$ in problem (\ref{formula:QPQC_Relaxed}) is neither an optimization variable  nor a parameter. Instead, it is an \emph{index} vector. In fact, the single constraint in problem (\ref{formula:QPQC_Relaxed}) is  a compact presentation for an \emph{infinite} number of constraints, where each constraint is  indexed by one choice of 
 $x\in\mathbb{R}^n$.   Second, if we set the scalars  $\varrho_{ij}$ and the vectors  $h_{m}$ to zero,  then problem (\ref{formula:QPQC_Relaxed}) reduces to the standard Lagrange dual problem associated with problem (\ref{formula:QPQC}), see \cite [Section 5.2]{Boyd}.  In that sense, problem (\ref{formula:QPQC_Relaxed}) can be seen as a \emph{generalized} dual problem for primal problem (\ref{formula:QPQC}), where the Lagrange multipliers corresponding to the linear inequality and linear equality constraints are affine rather than scalar \cite{Laurent_Book_2009}.  Third, the second line in (\ref{formula:QPQC_Relaxed}) involves multiplying every linear inequality constraint by itself and every other linear inequality constraint. Fourth,  the expression on the left hand side in the inequality constraints in (\ref{formula:QPQC_Relaxed}) is a quadratic function of index vector $x$.

Problem  (\ref{formula:QPQC_Relaxed}) is a relaxation of problem (\ref{formula:QPQC}), because any $\Lambda$ that satisfies the constraints in problem (\ref{formula:QPQC_Relaxed}) gives an \emph{upper bound} for the optimal objective value of the maximization in (\ref{formula:QPQC}).  In that sense,  problem (\ref{formula:QPQC_Relaxed}) seeks to find the lowest, i.e., the best, such upper bound \cite[Section 4.3]{SOS_Manual}.    The difference between the provided upper bound from (\ref{formula:QPQC_Relaxed}) and the true optimal objective value of problem (\ref{formula:QPQC}) is referred to as the \emph{relaxation gap}. In this paper, the relaxation gap is presented in percentage by dividing it by the true  optimal objective value of problem (\ref{formula:QPQC}).    If the resulting optimal $\Lambda$ is equal to the optimal objective value in (\ref{formula:QPQC}), then the relaxation is exact, and the \emph{relaxation gap} is zero. 
 For every  $x\in\mathbb{R}^n$ that is feasible in strategic bidding problem (\ref{formula:QPQC}), $X=xx^T$ is feasible in the proposed relaxation problem (\ref{formula:QPQC_Relaxed}). Thus,  the  infeasibility of problem (\ref{formula:QPQC_Relaxed}) is a certificate of infeasibility for problem (\ref{formula:QPQC}).
\color{black}


Problem (\ref{formula:QPQC_Relaxed}) is a convex optimization problem  because the objective function is linear and the feasible set is  convex.    However, since this problem has an \emph{infinite} number of constraints, i.e., one constraint for any $x \in \mathbb{R}^n$,    it is not a computationally tractable problem in its current form.  Therefore, next, we derive a tractable  representation for problem (\ref{formula:QPQC_Relaxed}).
\begin{lemma}\label{lemma:matrix}
Building upon the fourth property of problem (\ref{formula:QPQC_Relaxed}) mentioned earlier, its  constraint  can be reformulated as
\begin{equation}\label{formula:Compact_C}
{\begin{bmatrix} 1 \\ x \end{bmatrix}}^T
\Upsilon
\begin{bmatrix} 1 \\ x \end{bmatrix}\geq 0 \ \ \ \ \ \  \forall x\in \mathbb{R}^n,
\end{equation}
where
\begin{equation}\label{formula:Sch_SOS}
\begin{aligned}
\Upsilon \triangleq \!\!\!\!\! &&& \begin{bmatrix}  \Lambda  & -f^T \\   -f &-F  \end{bmatrix}
-\sum_{i}\alpha_i\begin{bmatrix} p_{i0} &  {p_i^T}/{2} \\    {p_i}/{2} & \mathbf{0} \end{bmatrix}- \\
%
%
 &&&\sum_{i=1}^I\sum_{j=1}^I\varrho_{ij}\begin{bmatrix} p_{i0}  \\   p_i\end{bmatrix}
 {\begin{bmatrix} p_{j0}  \\   p_j\end{bmatrix}}^T - \sum_{z=1}^Z\beta_z \begin{bmatrix}   0 & q_z^T \\ q_z & Q_z \end{bmatrix}-\\
 &&&\sum_{m=1}^Mh_{m0}\begin{bmatrix}1 \\ \mathbf{0}\end{bmatrix}{\begin{bmatrix} v_{m0}  \\  v_m\end{bmatrix}}^T-\sum_{m=1}^M\sum_{l=1}^{n} h_{ml}\begin{bmatrix}0 \\ e_l \end{bmatrix}{\begin{bmatrix} v_{m0}  \\  v_m\end{bmatrix}}^T.
\end{aligned}
\end{equation}
\end{lemma}
Here, $h_{ml}$ denotes the $l$th element of $h_m, \forall m$. 

From \cite[Excercise 3.32]{Thomas_2012}, a quadratic polynomial in $x$ such as the one on the left hand side of (\ref{formula:Compact_C}) in Lemma \ref{lemma:matrix} is always non-negative, if and only if it can be written as the sum of squares of some other polynomials \cite[Definition 3.24]{Thomas_2012}. From this, together with the analysis in   \cite[Section 3.1.4]{Thomas_2012},  the infinite number of  constraints in (\ref{formula:Compact_C}) is \emph{equivalent} to the following single matrix inequality constraint:
\begin{equation}\label{formula:Const_Semidef}
\Upsilon\succeq 0.
\end{equation}
 By replacing the  constraints in (\ref{formula:QPQC_Relaxed}) with  the one  in  (\ref{formula:Const_Semidef}), we express problem (\ref{formula:QPQC_Relaxed}) in  the following equivalent form:
 \begin{equation}\label{formula:SDP_Dual}
 \begin{aligned}
&\underset{\begin{subarray}{c}\Lambda,\alpha_i,\varrho_{ij}, \beta_z\\h_{ml},h_{m0} \end{subarray}}{\textbf{minimize} }       &&  \Lambda  \\
&\textbf{subject to} \!\!\! &&   \Upsilon\succeq 0.
\end{aligned}
 \end{equation}
 Problem (\ref{formula:SDP_Dual}) is a semidefinite program (SDP), which can be solved using  convex programming  tools such as Mosek \cite{Mosek}.

 \vspace{.1cm}
 \subsection{Reduced Computation Complexity}
 In this section, we reformulate  problem (\ref{formula:QPQC}) to significantly reduce the number of variables in problem  (\ref{formula:SDP_Dual}). This is done by systematically eliminating all linear equality constraints in problem (\ref{formula:QPQC}). First, we note that from
 \cite[pp. 46]{Longman},  set
\begin{equation}\label{formula:General_Linear}
\{x \ | \ v_m^Tx+v_{m0} =0,  \ \ \  \forall m\},
\end{equation}
is equivalent to set
\begin{equation}\label{formula:Short_Set}
\{Oy+\bar{x}\ | \ y\in\mathbb{R}^r\},
\end{equation}
 where
\begin{equation}\label{formula:Matlab_Command}
\begin{aligned}
& r \triangleq \text{Rank}\left( [v_1, \ldots, v_M]\right), \ \ O\triangleq \text{Null}({[v_1\cdots v_M]}^T) \\
& \: \ \ \ \ \ \ \ \ \ \bar{x}\triangleq  {[v_1\cdots v_M]}^T \texttt{\char`\\}{[v_{10}\cdots v_{M0}]}^T,
\end{aligned}
\end{equation}
 Here, the  matrix operator $\text{Null}( \cdot)$, which is also a command in Matlab \cite{Matlab_Linear}, returns an orthonormal basis for the null space of its argument matrix, obtained from its singular value decomposition.  Moreover, the operator ${\char`\\}$ which is also a command in Matlab \cite{Matlab_Linear}, returns an arbitrary member of the set (\ref{formula:General_Linear}).  From (\ref{formula:General_Linear}) and (\ref{formula:Short_Set}), we replace  optimization problem (\ref{formula:QPQC}) with the following equivalent optimization problem:
 \begin{equation}\label{formula:QPQC_Reduced}
\begin{aligned}
&\underset{y}{\textbf{maximize} }\ {\big(Oy+\bar{x}\big)}^TF\big(Oy+\bar{x}\big)+2f^T\big(Oy+\bar{x}\big)      && \\
&\textbf{subject to}  \\
                   & p_i^T\big(Oy+\bar{x}\big)+p_{i0}\geq 0      && \!\!\!\!\!\!\!\!\!\! \forall i      \\
&\big(Oy+\bar{x}\big)^TQ_z\big(Oy+\bar{x}\big)+2q_z^T\big(Oy+\bar{x}\big)=0    && \!\!\!\!\!\!\!\!\!\! \forall z.
\end{aligned}
\end{equation}
Note that, the above problem does \emph{not} have any linear equality constraint. While problem (\ref{formula:QPQC}) has $n$ variables, problem (\ref{formula:QPQC_Reduced}) has $r$ variables, where, in practice, $r \ll n$. Once we solve problem (\ref{formula:QPQC_Reduced}) and obtain its optimal solution $y^\star$, the optimal solution of problem (\ref{formula:QPQC}) is readily obtained as
\begin{equation}\label{formula:x_star}
{x^\star} = O {y^\star} + \bar{x}.
\end{equation}
 Similar to problem (\ref{formula:QPQC}), problem (\ref{formula:QPQC_Reduced}) is also a QCQP; therefore,  we can repeat the analysis in  Section \ref{sec:SolutionA} and introduce the following convex relaxation associated with problem (\ref{formula:QPQC_Reduced}):
 \begin{equation}\label{formula:SDP_Dual2}
 \begin{aligned}
&\underset{\begin{subarray}{c}\Lambda,\alpha_i,\varrho_{ij} \beta_z \end{subarray}}{\textbf{minimize} }       &&  \Lambda  \\
&\textbf{subject to} \!\!\! &&
{\Omega}^T  \Psi \
\Omega
\succeq 0,
\end{aligned}
 \end{equation}
 where
 \begin{equation}\label{formula:Psi}
\begin{aligned}
\Psi \triangleq &&& \begin{bmatrix}  \Lambda  & -f^T \\   -f &-F  \end{bmatrix}
-\sum_{i}\alpha_i\begin{bmatrix} p_{i0} &  {p_i^T}/{2} \\    {p_i}/{2} & \mathbf{0} \end{bmatrix}-\\
 &&&\sum_{i=1}^I\sum_{j=1}^I\varrho_{ij}\begin{bmatrix} p_{i0}  \\   p_i\end{bmatrix}
 {\begin{bmatrix} p_{j0}  \\   p_j\end{bmatrix}}^T - \sum_{z=1}^Z\beta_z \begin{bmatrix}   0 & q_z^T \\ q_z & Q_z \end{bmatrix},
\end{aligned}
\end{equation}
and
\begin{equation}
\Omega\triangleq \begin{bmatrix}
1 & \mathbf{0} \\ \bar{x} & O
\end{bmatrix}.
\end{equation}
Here, matrix $\Psi$ is a reduced version of matrix $\Upsilon$, where the optimization variables $h_{ml}$ and $h_{m0}$ are eliminated.  Similar to problem (\ref{formula:SDP_Dual}), problem (\ref{formula:SDP_Dual2}) is also an SDP. However, while problem (\ref{formula:SDP_Dual}) has ${n(n+1)/2}$ variables, problem (\ref{formula:SDP_Dual2}) has ${r(r+1)/2}$ variables. For example, for the case of the MPEC problem in Section \ref{sec:SingleTimeSlot}, the number of variables corresponding to problems (\ref{formula:SDP_Dual}) and (\ref{formula:SDP_Dual2}) are 11476 and 2016, respectively. This means 82\% drop in the number of variables.
 %
 %
 %

\vspace{.1cm}
\subsection{Recovery of  Original Optimization Variables}\label{sec:Recovery}
 In this section, we explain how we can recover a solution $y$ for problem (\ref{formula:QPQC_Reduced}) by solving its convex relaxation in (\ref{formula:SDP_Dual2}). A 
 solution $x$ for problem (\ref{formula:QPQC}) is then obtained 
 from $y$ using (\ref{formula:x_star}).  
 Suppose strong duality holds for the SDP in (\ref{formula:SDP_Dual2}), which is a convex optimization problem. Accordingly, problem (\ref{formula:SDP_Dual2}) and its dual problem, which itself is an SDP as shown below, have equal optimal objective values:
\begin{equation}\label{formula:SDP_Convert2}
\begin{aligned}
&\underset{ Y\in\mathbb{S}^{r+1}}{\textbf{maximize}}     \ \ \ tr \bigg(
{\Omega}^T
\begin{bmatrix}  0  & f^T \\   f &F  \end{bmatrix}
\Omega Y\bigg)  && \\
&\textbf{subject to } && \\
& Y_{11}=1  &&\\
& tr\bigg(
{\Omega}^T
\begin{bmatrix} p_{i0} &  p_i^T/2 \\    p_i/2 & \mathbf{0} \end{bmatrix}
\Omega Y
\bigg)\geq 0  && \forall i\\
  &tr\bigg(
  {\Omega}^T
  \begin{bmatrix} p_{i0}  \\   p_i\end{bmatrix}
 {\begin{bmatrix} p_{j0}  \\   p_j\end{bmatrix}}^T
    \Omega Y
 \bigg)\geq 0 &&    \forall i, j\\
  &tr\bigg(
  {\Omega}^T
  \begin{bmatrix}   0 & q_z^T \\ q_z & Q_z \end{bmatrix}
 \Omega Y\bigg)=0 &&  \forall z \\
 &{Y}\succeq 0. &&
\end{aligned}
\end{equation}
 Therefore, the above dual problem is still a convex relaxation of problem (\ref{formula:QPQC_Reduced}). Next, suppose  matrix ${Y^\star}$ denotes the optimal variable in problem (\ref{formula:SDP_Convert2}).
The following theorem explains the case where  the above convex relaxation  is \emph{exact}:
\vspace{.1cm}
\begin{theorem}\label{theorem:Optimality}
Suppose we obtain vector $y^\star \in \mathbb{R}^r$ from  matrix ${Y^\star}$  by taking the \emph{first column} of $Y^\star$ as follows:
\begin{equation}\label{formula:xstar}
\begin{bmatrix}  1 \\  {y^\star}  \end{bmatrix}=Y^\star e_1.
\end{equation}
 If $\text{Rank}(Y^\star) = 1$,  then  ${y^\star}$ is  the optimal solution of problem (\ref{formula:QPQC_Reduced}), and $x^\star$ in (\ref{formula:x_star}) is the optimal solution of problem (\ref{formula:QPQC}).
 \end{theorem}

\vspace{.1cm}
\noindent
The proof of Theorem \ref{theorem:Optimality} is given in  the Appendix. While Theorem \ref{theorem:Optimality} is promising, in practice, we often have $\text{Rank}(Y^\star) > 1$. Fortunately, even in that case, the approach in (\ref{formula:xstar}) gives a good \emph{approximate} solution for problem (\ref{formula:QPQC_Reduced}). That being said, there are still many  cases where such approximation is not  \emph{feasible}. Specially, ${y^\star}$ may not satisfy  all the quadratic equality  constraints in  (\ref{formula:QPQC_Reduced}).
 %
 %
Therefore, we need a mechanism to adjust  $y^\star$ from (\ref{formula:xstar}) to make it feasible.   Such mechanisms are often customized for  particular QCQP formulations, see \cite[Section IV-C]{Mohsenian_Miet_2011} for an example in Communications. In our case, we rather use the fact that the quadratic equality constraints in  (\ref{formula:QPQC_Reduced}) are all due to  \emph{complimentary slackness}, and hold the particular structure in (\ref{formula:Slack_Feature}). 
Accordingly, we propose Algorithm  \ref{algorithm:1} to
derive a feasible solution $y^\star$ from $Y^\star$.  The feasibility aspect  of solution from Algorithm \ref{algorithm:1} is analytically guaranteed, and its optimality is  shown to often be  exact  through  extensive numerical  Case Studies in Section \ref{sec:Case_Studies}.

 From the model in (\ref{formula:Slack_Feature}), the last constraint in (\ref{formula:QPQC_Reduced}) can be rewritten as
 \begin{equation}({d_z}^T{(O{y}+\bar{x})}+2)(q_z^T(O{y}+\bar{x}))=0  \ \ \ \forall z.
 \end{equation}
%
 %
 Therefore, we can express the last constraint in (\ref{formula:QPQC_Reduced}) as 
 \begin{equation}\label{formula:DCP}
{d_z}^T\big(O{y}+\bar{x}\big)+
2=0  \ \ \text{ or } \ \ q_z^T\big(O{y}+\bar{x}\big)=0, \ \ \; \forall z.
\end{equation}
Now, suppose for one quadratic equality constraint index $z$, neither  of the two equalities in (\ref{formula:DCP}) holds for $y = y^\star$, making $y^\star$ an infeasible solution to problem (\ref{formula:QPQC_Reduced}). But suppose there exists a small $\epsilon > 0$ and another number $\Delta \gg \epsilon$, for which
\begin{equation}\label{formula:condition1}
\begin{aligned}
| q_z^T\big(O{y^\star}+\bar{x}\big)| \leq \epsilon \ \ \text{ and } \ \ |{d_z}^T\big(O{y^\star}+\bar{x}\big)+
2| \geq \Delta.
 \end{aligned}
\end{equation}
 In that case, it is likely that at optimality we
  have
\begin{equation}\label{formula:Tem_Cons1}
q_z^T\big(O{y}+\bar{x}\big)=0.
\end{equation}
One can also make the  opposite argument. That is, if
\begin{equation}\label{formula:condition2}
|q_z^T\big(O{y^\star}+\bar{x}\big)| \geq \Delta \ \  \text{ and } \ \  |{d_z}^T\big(O{y^\star}+\bar{x}\big)+
2| \leq  \epsilon,
\end{equation}
then, it is likely that at optimality  we have
\begin{equation}\label{formula:Tem_Cons2}
{d_z}^T\big(O{y}+\bar{x}\big)+
2=0.
\end{equation}
Therefore, if it turns out that (\ref{formula:condition1}) holds for a specific index $z$, then we can  replace the corresponding   complimentary slackness constraint in (\ref{formula:QPQC_Reduced})  which is \emph{non-convex}, with its equivalent-at-optimality \emph{linear}  constraint in (\ref{formula:Tem_Cons1}).  Similarly, if (\ref{formula:condition2}) holds for a specific $z$, the corresponding complimentary
 slackness constraint in problem (\ref{formula:QPQC_Reduced}) is replaced with  (\ref{formula:Tem_Cons2}).%
%
%

The above argument is the foundation of  Algorithm \ref{algorithm:1}.   Once we encounter an infeasible solution $y^\star$ in Line 3, we first initialize the values of parameters $\Delta$ and $\epsilon$ in Line 4, and then we go through \emph{iterations} of augmenting problem (\ref{formula:QPQC_Reduced}) in Lines 5 to 11 until we obtain a feasible solution. In the first iteration, we deal with a version of problem (\ref{formula:QPQC_Reduced}) in which we have removed several  complimentary slackness constraints through Lines 5 to 9. Therefore, solving the \emph{MILP-equivalent} of such  augmented  problem in Line 10 is a light  task. Next, as  we keep iterating through Lines 5 to 11, we decrease $\epsilon$, and we choose to keep more original complimentary slackness constraints in problem (\ref{formula:QPQC_Reduced}), until the augmented problem (\ref{formula:QPQC_Reduced}) becomes feasible. Accordingly, the computation time in solving the MILP-equivalent of problem (\ref{formula:QPQC}) will gradually grow as we iterate.  However, 
as we will see in Section \ref{sec:SingleTimeSlot}, in practice, we often need to iterate very few times; therefore, in general, the computation time for Algorithm \ref{algorithm:1} is much lower compared to the standard MILP approach in \cite{ruiz2009pool,MILP_Example_1,MILP_Example_2}.

 In summary,   Algorithm \ref{algorithm:1} exploits the solution that comes  from the proposed relaxation problem in (\ref{formula:SDP_Convert2}) in order to reduce the computation time in solving problem (\ref{formula:QPQC}).   The solution of Algorithm \ref{algorithm:1} is guaranteed to be feasible   to problem (\ref{formula:QPQC}), due to Steps 3 and 12 in Algorithm 1. However,  neither  the  computation time  nor the optimality of Algorithm \ref{algorithm:1} is guaranteed. Nevertheless, the numerical examples in Section \ref{sec:Case_Studies} suggest  that Algorithm \ref{algorithm:1} often  performs  very effectively in solving problem (\ref{formula:QPQC}), with high optimality and low computation time.    As for the convergence of Algorithm \ref{algorithm:1}, we note that,  it iteratively solves a finite number of MILPs one-after-one until one does  converge.  In the worst case scenario,  Algorithm \ref{algorithm:1} would end up solving the original    MILP reformulation of problem (\ref{formula:QPQC}) based on \cite{ruiz2009pool}, which is guaranteed to converge to a feasible solution,  but it may take a long time to do so.  This is because problem (\ref{formula:QPQC}) is a reformulation of problem (\ref{formula:MPEC_KKT}), and by construction problem (\ref{formula:MPEC_KKT}) is  always
feasible.

\begin{algorithm}[t]
\caption{}
\label{algorithm:1}
\begin{algorithmic}
\STATE \hspace{.02cm} 1: Solve convex relaxation problem (\ref{formula:SDP_Convert2}) and obtain $Y^\star$.
\STATE \hspace{.02cm} 2: Obtain $y^\star$  from $Y^\star$ using (\ref{formula:xstar}).
\STATE \hspace{.02cm} 3: \textbf{if} $y^\star$ is feasible to problem (\ref{formula:QPQC_Reduced}) \textbf{then} exit.
 \STATE \hspace{.02cm} 4: Set $\Delta=1$ and $\epsilon=0.1$.
\STATE \hspace{.02cm} 5: \textbf{for} each complimentary slackness   constraint $z$  \textbf{do}
\STATE \hspace{.02cm} 6:  \hspace{.7cm} \textbf{if}   condition (\ref{formula:condition1}) holds for $y = y^\star$ \textbf{then}
\STATE \hspace{.02cm} 7:  \hspace{1.05cm} Replace  constraint  $z$ in (\ref{formula:QPQC_Reduced}) with (\ref{formula:Tem_Cons1}).
\STATE\hspace{.02cm}  8:  \hspace{.7cm} \textbf{if}   condition (\ref{formula:condition2}) holds for $y = y^\star$ \textbf{then}
\STATE  \hspace{.02cm} 9:  \hspace{1.05cm} Replace  constraint $z$ in (\ref{formula:QPQC_Reduced}) with (\ref{formula:Tem_Cons2}).
\STATE  10: Solve the MILP equivalent of problem (\ref{formula:QPQC_Reduced}), see \cite{ruiz2009pool}.
\STATE  11: Set $\epsilon=\epsilon-0.01.$
\STATE 12: \textbf{if} the MILP equivalent  is infeasible \textbf{then} Go to Step 5.
\end{algorithmic}
\end{algorithm}

\vspace{.11cm}
\section{Multiple Time Slots and Random Scenarios}\label{sec:Multi}

In practice,    problem (\ref{formula:MPEC_KKT}) may need to be solved over  $T\geq 1$ time slots,  e.g., over 24 hourly time slots in a day-ahead market. Also, one may often need to address
 uncertainty by taking into account $K \geq 1$ random scenarios. In that case,  the price and energy bid  parameters of generators and loads and also all the variables in  MPEC problem  are indexed by $t$ and $k$.  For example, $x_k[t]$ means the vector of the original optimization variables $x$ indexed at time slot $t$ and random scenario $k$. Hence, we can extend the MPEC problem formulation in (\ref{formula:MPEC_KKT}) and present it in vector-format as \cite{ruiz2009pool}:
\begin{equation}\label{formula:QPQC_Multi}
\begin{aligned}
&\underset{x_{k}[t]}{\textbf{maximize} } \ \  \sum_{t=1}^{T}\sum_{k=1}^{K}     {x_{k}[t]}^T \frac{F_k[t]}{K}x_{k}[t]+\sum_{t=1}^{T}\sum_{k=1}^{K}2\frac{{f_k[t]}^T}{K}x_{k}[t]  &&\\
& \textbf{subject to}  \\
 & {{p_{i,k}[t]}}^Tx_{k}[t]+{{p_{i0,k}[t]}}\geq 0 &&\!\!\!\!\!\!\!\!\!\!\!\!\!\!\!\!\!\!\!\!\!\!\!\!\!\!\!\!\!\!\!\!\!\!\!\!\!\!\!\!\!\!\!\!\!\! \forall t,  k, i\\
&{{{v_{m,k}[t]}}}^Tx_{k}[t]+{{v_{m0,k}[t]}}=0 &&\!\!\!\!\!\!\!\!\!\!\!\!\!\!\!\!\!\!\!\!\!\!\!\!\!\!\!\!\!\!\!\!\!\!\!\!\!\!\!\!\!\!\!\!\!\!
\forall t,  k,   m\\
&  {x_{k}[t]}^T{{Q_{z,k}[t]}}x_{k}[t]+2{{{q_{z,k}[t]}}}^Tx_{k}[t]=0 &&  \!\!\!\!\!\!\!\!\!\!\!\!\!\!\!\!\!\!\!\!\!\!\!\!\!\!\!\!\!\!\!\!\!\!\!\!\!\!\!\!\!\!\!\!\!\!
\forall t,  k, z \\
&e_l^Tx_{k}[t]-e_l^Tx_{k}[t-1]+\Gamma\geq0   && \!\!\!\!\!\!\!\!\!\!\!\!\!\!\!\!\!\!\!\!\!\!\!\!\!\!\!\!\!\!\!\!\!\!\!\!\!\!\!\!\!\!\!\!\!\!   \forall t, k,  \exists l\\
&e_l^Tx_{k}[t-1]-e_l^Tx_{k}[t]+\Gamma\geq0   && \!\!\!\!\!\!\!\!\!\!\!\!\!\!\!\!\!\!\!\!\!\!\!\!\!\!\!\!\!\!\!\!\!\!\!\!\!\!\!\!\!\!\!\!\!\!   \forall t, k,  \exists l \\
& e_l^Tx_{k}[t]-e_l^Tx_{1}[t]=0
     && \!\!\!\!\!\!\!\!\!\!\!\!\!\!\!\!\!\!\!\!\!\!\!\!\!\!\!\!\!\!\!\!\!\!\!\!\!\!\!\!\!\!\!\!\!\!  \forall t,   k,  \exists l.
\end{aligned}
\end{equation}
 The notation $\forall t, k$  in the constraints of problem (\ref{formula:QPQC_Multi}) indicates that the corresponding constraints hold for all the time slots and  all the scenarios within their corresponding  ranges, i.e., $t=1,\cdots\!,T$ and $k=1,\cdots\!,K$. Also, the notation $\exists l$ indicates that the  constraint  holds only for strategic generators.  \color{black}
The first three constraints in (\ref{formula:QPQC_Multi}) are simply the extensions of the  constraints in problem (\ref{formula:QPQC}),  across time slots and random scenarios. The fourth and fifth  constraints in  (\ref{formula:QPQC_Multi})   includes the   ramp constraints for strategic generators, where in each case the index $l$ and accordingly the basis $e_l$ are selected such that $e_l^T x_k[t]$ indicates the generation output of a particular strategic generator at time slot $t$ and random scenario $k$. Finally, the sixth constraint in (\ref{formula:QPQC_Multi}) is used to make sure that the bids of the strategic generators are the same across all random scenarios, where in each case the index $l$ and accordingly the basis $e_l$ are selected such that $e_l^T x_k[t]$ indicates the price bid of a particular strategic generator at time slot $t$ and random scenario $k$. 

\vspace{.11cm}
\subsection{Immediate Solution Approach}\label{sec:Immediate}
 Just like problem (\ref{formula:QPQC}), problem  (\ref{formula:QPQC_Multi}) is a QCQP. However, the size of the optimization vector in (\ref{formula:QPQC_Multi}) is
 $TK$ times the size of the optimization vector in problem (\ref{formula:QPQC}). One approach to solve problem (\ref{formula:QPQC_Multi}) is to follow exactly the same analysis in Section \ref{sec:Solution}.  This is done by expanding the inequality constraint in (\ref{formula:QPQC_Relaxed}) to also include the last three constraints in problem (\ref{formula:QPQC_Multi}). Specifically, since the last three constraints in (\ref{formula:QPQC_Multi}) are linear, their corresponding Lagrange multipliers in (\ref{formula:QPQC_Relaxed}) would be affine, just like the case of the linear constraints in problem (\ref{formula:QPQC}), please refer to the second and the third properties of problem (\ref{formula:QPQC_Relaxed}) that we discussed in Section \ref{sec:SolutionA}.

Once problem (\ref{formula:QPQC_Relaxed}) is updated as we explained above,  we would then follow  the rest of the analysis in Section \ref{sec:Solution}
 and end up with solving an
 SDP similar to the one in  (\ref{formula:SDP_Convert2}).  While in this approach we would achieve a  convex relaxation for problem (\ref{formula:QPQC_Multi}), the matrix domain of the resulting  SDP problem  would  be   $\mathbb{S}^{TK r+1}$, which means having   $TK r(TK r+1)/2$  scalar variables. Unfortunately, the  number of constraints in such SDP grows in proportional   to  $T^2K^2$. In other words,  even though the problem itself remains convex, its size will grow exponentially as the number of time slots and random scenarios grows. As a result, such convex relaxation may impose
 huge computation burden and may  not be practical.

\vspace{.11cm}
\subsection{Alternative Solution Approach}\label{sec:Alternative}

 In this section, we propose an alternative  convex  relaxation approach  to solve problem (\ref{formula:QPQC_Multi})  to tackle the \emph{curse of dimensionality} in  the number of time slots and random scenarios.    Again, we start  by expanding the inequality constraint in (\ref{formula:QPQC_Relaxed}) to also include the last three constraints in problem (\ref{formula:QPQC_Multi}). However, as opposed to the approach in Section \ref{sec:SolutionA}, where we would use affine Lagrange multipliers for these three new sets of linear constraints, we would use only scalar Lagrange multipliers, just like in the standard Lagrange dual problem formulation \cite[Section 5.2]{Boyd}. This would, in presence of large $T$ and $K$, significantly reduce the number of additional Lagrange multipliers in the extension of problem (\ref{formula:QPQC_Relaxed}); and accordingly the number of variables in problem (\ref{formula:Multi}). The rest of the analysis would be similar to Section \ref{sec:Solution}. Here, we only show the final convex relaxation problem that we must solve:
%
%
%
%
%
%
%
\begin{equation}\label{formula:Multi} \nonumber
\begin{aligned}
& \underset{ y_k[t],Y_k[t]}{\textbf{maximize}}  \    \frac{1}{K}\sum_{t=1}^{T}\sum_{k=1}^{K} tr \bigg(
{{\Omega_k[t]}}^T
\begin{bmatrix}  0  & {f_{k}[t]}^T \\   f_{k}[t] &F_{k}[t]  \end{bmatrix}
{\Omega_k[t]} Y_k[t]\bigg)   &&\!\!\!\!\!\!\!\!\!\!\!\!\!\!\!\!\!\!\!\!\!\!\!\!\!\!\!\!\!  \\
&\textbf{subject to } &&  \\
%
%
& Y_{11,k}[t]=1  &&\!\!\!\!\!\!\!\!\!\!\!\!\!\!\!\!\!\!\!\!\!\!\!\!\!\!\!\!\!\!\!\! \forall t,  k\\
& tr\bigg({{\Omega_k[t]}}^T
\begin{bmatrix} {{p_{i0,k}[t]}} &  \frac{{p_{i,k}[t]}^T}{2}\\    \frac{{p_{i,k}[t]}}{2} & \mathbf{0} \end{bmatrix}
{\Omega_k[t]}Y_k[t]
\bigg)\geq 0  &&\!\!\!\!\!\!\!\!\!\!\!\!\!\!\!\!\!\!\!\!\!\!\!\!\!\!\!\!\! \!\!\! \forall t, k ,  i \\
  &tr\bigg( {{\Omega_k[t]}}^T \begin{bmatrix} {{p_{i0,k}[t]}}  \\   {{p_{i,k}[t]}}\end{bmatrix}
 {\begin{bmatrix} p_{j0,k}[t]  \\   p_{j,k}[t]\end{bmatrix}}^T
  \!\! {\Omega_k[t]} Y_k[t]
 \bigg)\geq 0, &&\!\!\! \!\!\!\!\!\!\!\!\!\!\!\!\!\!\!\!\!\!\!\!\!\!\!\!\!\!\!\!\! \forall t, k,  i,  j \\
\end{aligned}
\end{equation}
\begin{equation}\label{formula:Multi}
\begin{aligned}
  &tr\bigg({{\Omega_k[t]}}^T
  \begin{bmatrix}   0 & {{{q_{z,k}[t]}}}^T \\ {{q_{z,k}[t]}} & {{Q_{z,k}[t]}} \end{bmatrix}
  {\Omega_k[t]}Y_k[t]
  \bigg)=0, &&  \!\!\! \forall t,  k,  z \\
 &Y_k[t]\succeq 0, &&  \!\!\! \forall t, k, \\
&\begin{bmatrix}  1 \\  {y_k[t]}  \end{bmatrix}=Y_k[t] e_1, &&  \!\!\! \forall t, k  \\
 &e_l^T{O_k[t]}^T\left(y_k[t]-y_k[t-1]\right)+\Gamma\geq0,    &&  \!\!\! \forall t, k,  \exists l \\
  &e_l^T{O_k[t]}^T\left(y_k[t-1]-y_k[t]\right)+\Gamma\geq0,   &&  \!\!\! \forall t, k,  \exists l \\
  &e_l^T{O_k[t]}^T\left(y_k[t]-y_1[t]\right)=0,  && \!\!\! \forall t, k,  \exists l, \ 
\end{aligned}
\end{equation}
where
\begin{equation}
{\Omega_k[t]}\triangleq \begin{bmatrix} 1 & \mathbf{0}\\ {\bar{x}_k[t]}  & {O_k[t]} \end{bmatrix}  \ \ \  \forall t, k.
\end{equation}

\begin{figure}[t]
\captionsetup{font=MyFigureFont}
\begin{center}
\scalebox{0.51}{\includegraphics*{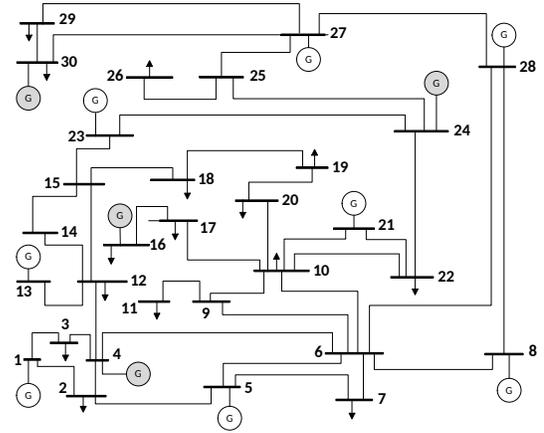}}\end{center}
\vspace{-0.2cm}\caption{The IEEE 30-bus test system  that we considered in our case studies. The generators in the strategic
generation firm are highlighted in gray.} \vspace{-0.3cm} \label{fig:IEEE30bus}
\end{figure}

\begin{table*}[t]
\caption{Load Data}
\begin{center}
\begin{tabular}{| c | c | c | c | c | c | c | c | c | c | c | c | c | c | c | c | c | c | c | c | c | c | c | c | c |}
\hline
\!\!\!\! & \multicolumn{24}{ c| }{Hourly Price Bids} \\ \cline{2-25}
\!\!Bus \#\!\! & $1$ &  $2$  &  $3$  &  $4$  &  $5$  &  $6$  &  $7$  &  $8$  &  $9$  &  $10$  &  $11$  &  $12$  &  $13$  &  $14$  &  $15$  &  $16$  &  $17$  &  $18$  &  $19$  &  $20$  &  $21$  &  $22$  &  $23$  &  $24$  \\ \hline
\! 26 \!\! & \!43.5\! & \!41.6\! & \!33.7\! & \!36.1\! & \!35.5\! & \!43.9\! & \!48.2\! & \!58.0\! & \!41.0\! & \!46.2\! & \!41.9\! & \!43.8\! & \!43.9 \! & \!45.0\! & \!44.0\! & \!42.5\! & \!48.4\! & \!58.4\! & \!63.0\! & \!72.4\! & \!65.7\! & \!59.1\! & \!52.7\! & \!48.7 \! \\ \hline
\! 29 \!\! & \!42.4\! & \!38.0\! & \!35.8\! & \!38.0\! & \!38.2\! & \!40.5\! & \!54.3\! & \!60.0\! & \!53.1\! & \!47.0\! & \!44.5\! & \!45.8\! & \!41.6\! & \!41.7\! & \!44.9\! & \!48.9\! & \!48.8\! & \!59.2\! & \!62.1\! & \!68.2\! & \!64.0\! & \!62.4\! & \!53.1\! & \!45.0 \! \\ \hline
\end{tabular}
\end{center}
\label{table:Load_Data}
\end{table*}

\begin{table*}[t]
\caption[caption]{Scaling Factors for Construction of Random Scenarios}
\begin{center}
\begin{tabular}{| c | c | c | c | c | c | c | c | c | c | c | c | c | c | c | c | c | c | c | c | c | c | c | c | c |}
\hline
 \multicolumn{20}{ |c| }{Scenario \#} \\ \cline{1-20}
 $1$  &  $2$   &  $3$  &  $4$  &  $5$  &  $6$  &  $7$  &  $8$   &  $9$   &  $10$  &  $11$  &  $12$  &  $13$  &  $14$ & $15$  &  $16$  &  $17$  &  $18$  &  $19$  &  $20$    \\ \hline
  \! 2.0 \! & \!1.9\! & \!1.8\! & \!1.7\! & \!1.6\! & \!1.5\! & \!1.4\! & \!1.3\! & \!1.2\! & \!1.1\! & \!1.0\! & \!0.9\! & \!0.8\! & \!0.7\! & \!0.6\! & \!0.5\! & \!0.4\! & \!0.3\! & \!0.2\! & \!0.1\!  \\ \hline
\end{tabular}
\end{center}
\label{table:Factors}
\end{table*}

\newcolumntype{C}[1]{>{\centering\let\newline\\\arraybackslash\hspace{0pt}}m{#1}}

\begin{table*}[htbp]
  \small
  \centering
  \caption{Computation Time  for Different  Separate Time Intervals \color{black} (Minutes)}
  \subtable[Proposed Approach]{
    \begin{tabular}{| c | C{.7cm} | C{.7cm} | C{.7cm} | C{.7cm} | C{.7cm} | C{.7cm} |   }
\hline
 & \multicolumn{6}{ c| }{Time Interval} \\ \cline{2-7}

\!\! K\!\! & $1$ &  $2$  &  $3$  &  $4$  &  $5$  &  $6$     \\ \hline

  \!\! 6   \!\! &    11   \!\! &  \!\!   15      \!\!  &   \!\!  11     \!\!  &  \!\!  12      \!\!  & \!\!  9      \!\!  &  10        \!\! \\ \hline
  \!\! 9  \!\!  &  13    \!\! &  \!\!   17    \!\!  &   \!\!  9     \!\!  &   \!\!  7     \!\!  & \!\!  13      \!\!   &  7       \!\! \\ \hline
  \!\! 8  \!\!  &  14  \!\! &  \!\!   19    \!\!  &  \!\!  18    \!\!  &   \!\!  15     \!\!  & \!\!  14       \!\!  &  11       \!\! \\ \hline
  \!\! 9  \!\!  &  11    \!\! &  \!\!   21   \!\!  &  \!\!  20     \!\!  &    \!\!  11     \!\!  & \!\! 17       \!\!  &  10       \!\! \\ \hline
  \!\! 10  \!\!  & 19  \!\! &  \!\!   21   \!\!  &  \!\!  12    \!\!  &    \!\!  10    \!\!   & \!\! 12       \!\!  &16       \!\!  \\ \hline   

\end{tabular}
  } \ \ \ \ \ \ \ \ \ \ \  
  \subtable[MILP Approach in \cite{ruiz2009pool}]{
    \begin{tabular}{| c | C{.7cm} | C{.7cm} | C{.7cm} | C{.7cm} | C{.7cm} | C{.7cm} | C{.7cm}   }
\hline
 & \multicolumn{6}{ c| }{ Time Interval} \\ \cline{2-7}

\!\! K\!\! & $1$ &  $2$  &  $3$  &  $4$  &  $5$  &  $6$     \\ \hline

  \!\! 6   \!\! &    236   \!\! &  \!\!   97      \!\!  &   \!\!  516     \!\!  &  \!\!  174      \!\!  & \!\!  35      \!\!  &  21        \!\! \\ \hline
  \!\! 7  \!\!  &  949    \!\! &  \!\!   119    \!\!  &   \!\!  287     \!\!  &   \!\!  285     \!\!  & \!\!  54      \!\!   &  57       \!\! \\ \hline
  \!\! 8  \!\!  &  2126  \!\! &  \!\!   537    \!\!  &  \!\!  1135    \!\!  &   \!\!  137     \!\!  & \!\!  23       \!\!  &  40       \!\! \\ \hline
  \!\! 9  \!\!  &  553    \!\! &  \!\!   -   \!\!  &  \!\!  329     \!\!  &    \!\!  470     \!\!  & \!\! 71       \!\!  &  70       \!\! \\ \hline
  \!\! 10  \!\!  & 2197  \!\! &  \!\!   -   \!\!  &  \!\!  -    \!\!  &    \!\!  587    \!\!   & \!\! 70       \!\!  &161       \!\!  \\ \hline   

\end{tabular}
  }
  \label{table:CMT}
\end{table*}

 Next, we highlight some of the key properties of problem (\ref{formula:Multi}). First, if $T = K = 1$, then problem (\ref{formula:Multi}) reduces to problem (\ref{formula:SDP_Convert2}), where the last three sets of constraints in (\ref{formula:Multi}) will disappear and the sixth constraint in (\ref{formula:Multi}) reduces to (\ref{formula:xstar}) in Theorem \ref{theorem:Optimality}. Second, the SDP problem in (\ref{formula:Multi}) has  a mix of matrix variables $Y_k[t]$ and vector variables $y_k[t]$. Third, the  number of  variables in  problem (\ref{formula:Multi}) is only $TKr(r+1)/2$, which grows only \emph{linearly} with respect to  either  the number of time slots $T$ or  the number of random scenarios $K$. As we will see in Sections \ref{sec:SingleScenario} and \ref{sec:SingleTimeSlot},  this  latter property plays a drastic role in lowering the computation time in our proposed approach, compared to the standard MILP approach in \cite{ruiz2009pool, MILP_Example_1,MILP_Example_2}.   Fourth,  matrices $Y_{k}[t]$ $\forall t,k$ are dense, i.e., not sparse. Therefore, the matrix completion methods such as the one in \cite{Vanden2015,Jabr2012} are not applicable to problem (\ref{formula:Multi}).  \color{black}

 As in Theorem \ref{theorem:Optimality}, if  $\text{Rank}(Y_k^\star[t])=1$,  $\forall t,k$, then the convex relaxation  in problem (\ref{formula:Multi}) is \emph{exact}, i.e., the optimal solutions of the original MPEC  problem in (\ref{formula:QPQC_Multi}) are obtained as
\begin{equation}
x^\star_k[t] = O y_k^\star[t] + \bar{x}_k[t],
\end{equation}
where $y^\star_k[t]$,  $\forall t, k$ is the optimal solution of problem (\ref{formula:Multi}). Again, in practice, $\text{Rank}(Y_k^\star[t])  >1$,  for several time slot $t$ and random scenario $k$ instances. In such cases, we can still use Algorithm \ref{algorithm:1},  where we replace Lines 1 and 2 with ``Solve Problem  (\ref{formula:Multi}) and obtain $Y^\star_k[t]$ and $y^\star_k[t]$ for all $t$ and $k$.''

\vspace{.11cm}
\section{Case Studies}\label{sec:Case_Studies}

\vspace{.1cm}
\subsection{Simulation Setup}
In this section, we assess the performance of the proposed approach  based on  the extended  IEEE 30 bus test system in \cite{Mohsenian_MPEC_2015}, see Fig. \ref{fig:IEEE30bus}, where the four generators in the strategic generation firm are highlighted using color gray. Here, the network includes 30 buses and
41 transmission lines. We have: 
$\mathcal{S}$ = \{4, 16, 24, 30\}. The transmission lines data,  generation data, and load energy bids data are the same as those  in  Tables I to  III in    \cite{Mohsenian_MPEC_2015}.  Specifically, the transmission line between bus \#2 and bus \#4 has a limited capacity of 0.2.  Each strategic generation unit has 1 GW capacity and the ramp parameter is $\Gamma=0.3$. The  cost vector for strategic generators is $c^T$= [45.84 47.84   55.56 63.88] \$/MWh.  All loads, except for those at buses 26 and 29,  submit a price bid of 72 \$/MWh for all 24 market operation hours. The hourly price bids  of the load  at bus 26 and bus 29 are as in Table \ref{table:Load_Data}.    As in \cite{ruiz2009pool}, we construct  20 random scenarios by  scaling the price bids of loads and non-strategic generators by using the 20 scaling factors  that are given  in Table \ref{table:Factors}.    All problems are solved using a single  Intel Xeon E5-2450-v2  CPU.      

 Problem (\ref{formula:Multi}) is solved using Yalmip \cite{Yalmip}, where Mosek \cite{Mosek} is the SDP solver. All MILP formulations are    solved using  Gurobi \cite{Gurobi}.   In all  case studies, the number of time slots  $T$ and the number of  random scenarios $K$ are selected such that, the MILP approach in \cite{ruiz2009pool} can converge in a  timely manner to allow us assess the optimality of our own design.   The optimality of our proposed approach   is measured based on the profit that is gained by  the generation firm,  after bidding the solution that comes from our proposed approach. In this regards,      
 for any solution  $x^\star$ that is feasible to the constraints of problem (\ref{formula:MPEC_KKT}), the optimality of $x^\star$ is defined as the numerical  value of the objective function in (\ref{formula:MPEC_KKT})  at  $x=x^\star$, divided by the true  optimal objective value of problem (\ref{formula:MPEC_KKT}).
 \color{black}

\vspace{.11cm}
\subsection{Impact of Increasing the Number of  Random Scenarios}\label{sec:SingleTimeSlot}
\begin{figure}[t]
    \centering \subfigure[]{  {\scalebox{0.45} {\includegraphics{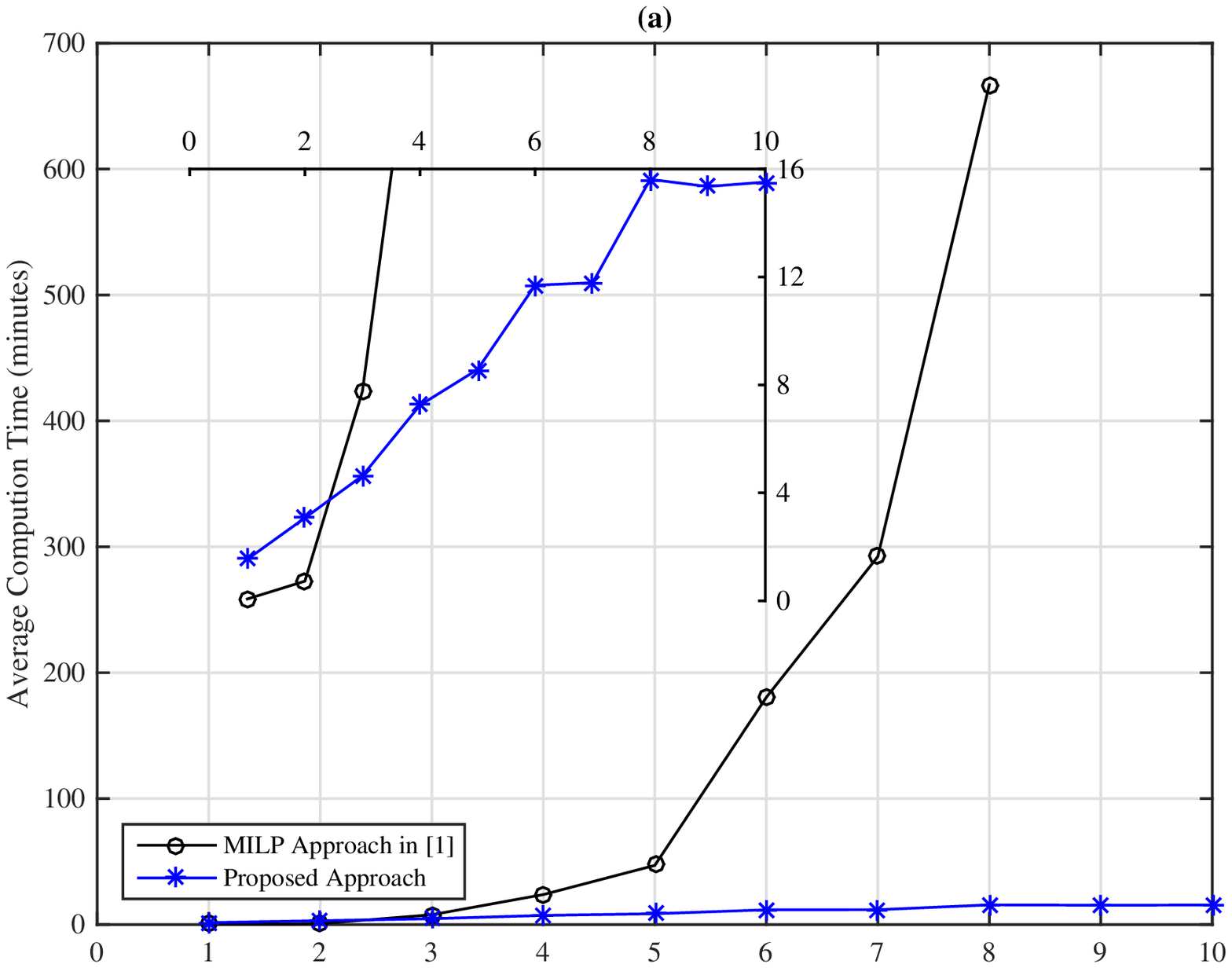}}} }
    \
      
 \vspace{-1.2cm}
 
    \centering \subfigure[]{  {\scalebox{0.45} {\includegraphics{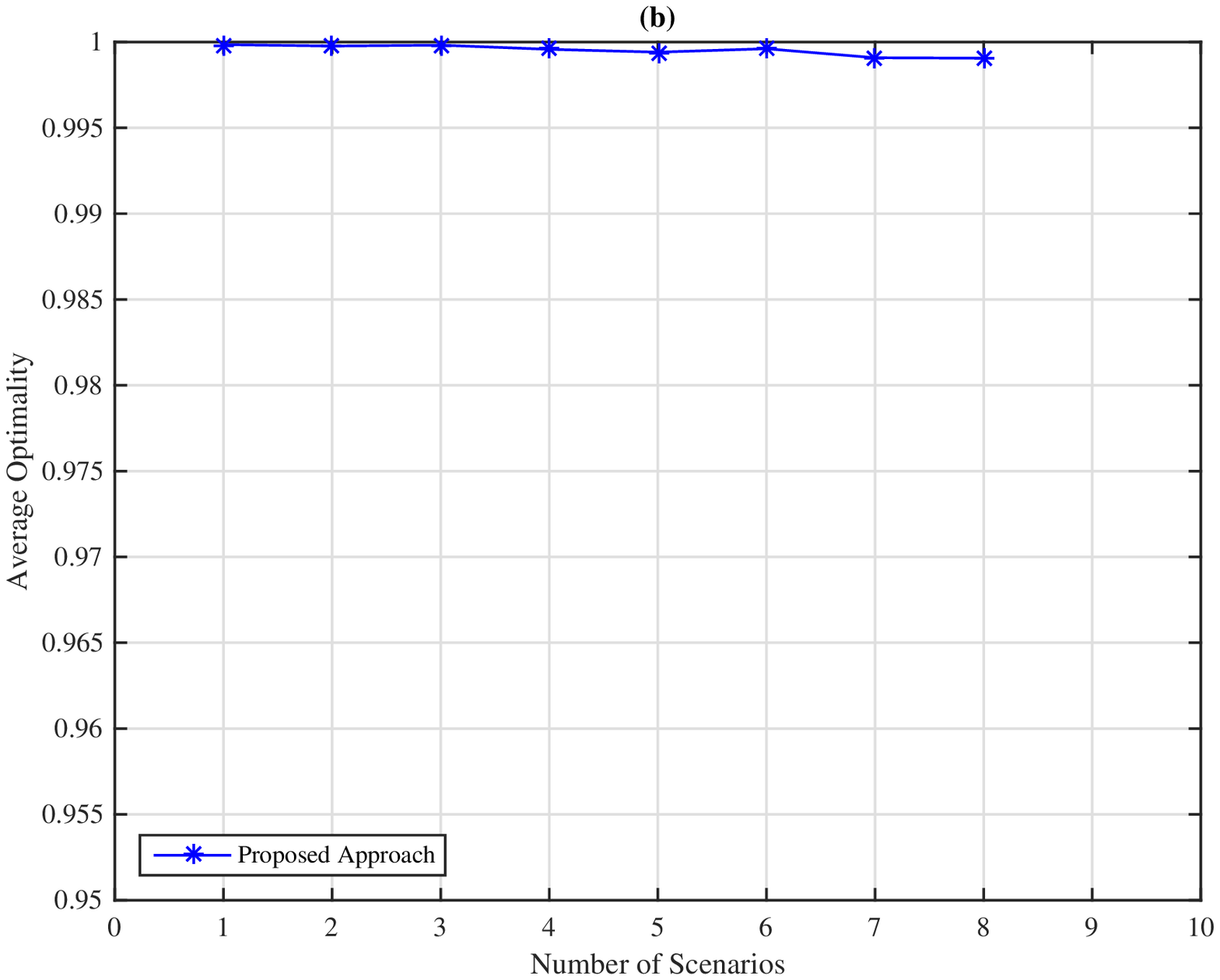}}} }
    \
      
    \vspace{-.5cm}
      \caption{The impact of increasing the number of random scenarios on the performance of the proposed approach and the MILP approach in \cite{ruiz2009pool}: (a) the computation time; (b) the optimality.}
      \vspace{-.0cm} 
 \label{fig:1}
\end{figure}

 Suppose  $T = 4$.    Fig. \ref{fig:1}(a) shows the average  computation time  versus the number of random scenarios $K$ for our  approach as well as for the  MILP approach in \cite{ruiz2009pool}.   Here, the average is taken across   six  MPEC problems for  six different time intervals of length four hours.   From Fig. \ref{fig:1}(a), we can see that as $K$ increases,  the computation time  for MILP approach in \cite{ruiz2009pool} grows exponentially while for our  proposed approach grows rather linearly. The difference between the two approaches becomes particularly significant where there are $K=6$ or more random scenarios.   For this range of random scenarios, the computation times are shown in Table \ref{table:CMT}. We can see that, when $K=10$, the MILP approach in \cite{ruiz2009pool} does not converge  for the second and third time intervals, even after  running for three days. In contrast,  our approach always converged in less than 21 minutes.  Interestingly,  the optimality of the solution that comes from our proposed approach  is always $96\%$ or better, for all the cases that are studied in this section.   For  example, where  $K = 8$,  the average computation time for  the approach in \cite{ruiz2009pool} and our approach are  667  minutes versus only about   16  minutes, respectively. This suggests an improvement factor   over 40.   Note that, we did  not  go beyond   $K = 10$  scenarios, mainly because the MILP approach in \cite{ruiz2009pool} could not converge in a timely manner for the larger   number of scenarios.  In particular, the MILP approach in \cite{ruiz2009pool} could not converge even after running the MILP algorithm  for three days.   Otherwise, as far as our proposed approach is concerned, we can handle larger $K$ in this case, if needed.

Next, we take a closer look at how Algorithm \ref{algorithm:1} behaves.  Out  of the $10\times6=60$ total case instances  that are analyzed in the case studies in  this Section, in 36 cases, the inner loop of Algorithm \ref{algorithm:1} was executed only once. In  24 cases, the inner loop of Algorithm \ref{algorithm:1} was iteratively executed between two to nine  times.
That being said,  Algorithm \ref{algorithm:1} never  iterated more than nine times between Step 5 and Step 12, and   never  ended up solving  the original problem in (\ref{formula:QPQC}) using the MILP approach \cite{ruiz2009pool}.  Of course, this may change in other test cases.

  %

\vspace{.11cm}
\subsection{Impact of Increasing the Scheduling Horizon}\label{sec:SingleScenario}

\begin{figure}[t]
    \centering \subfigure[]{  {\scalebox{0.45} {\includegraphics{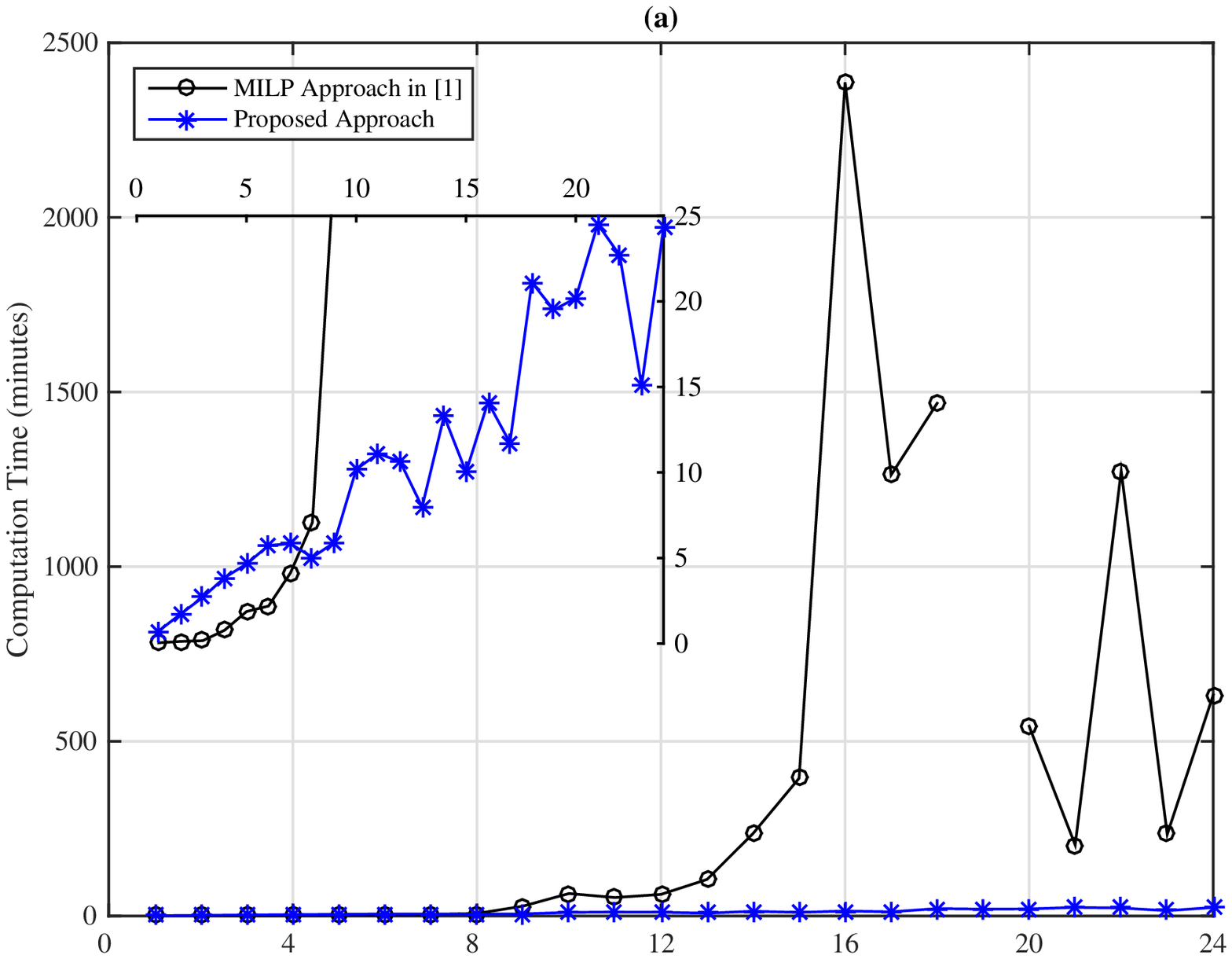}}} }
    \
      
 \vspace{-1.2cm}
 
    \centering \subfigure[]{  {\scalebox{0.45} {\includegraphics{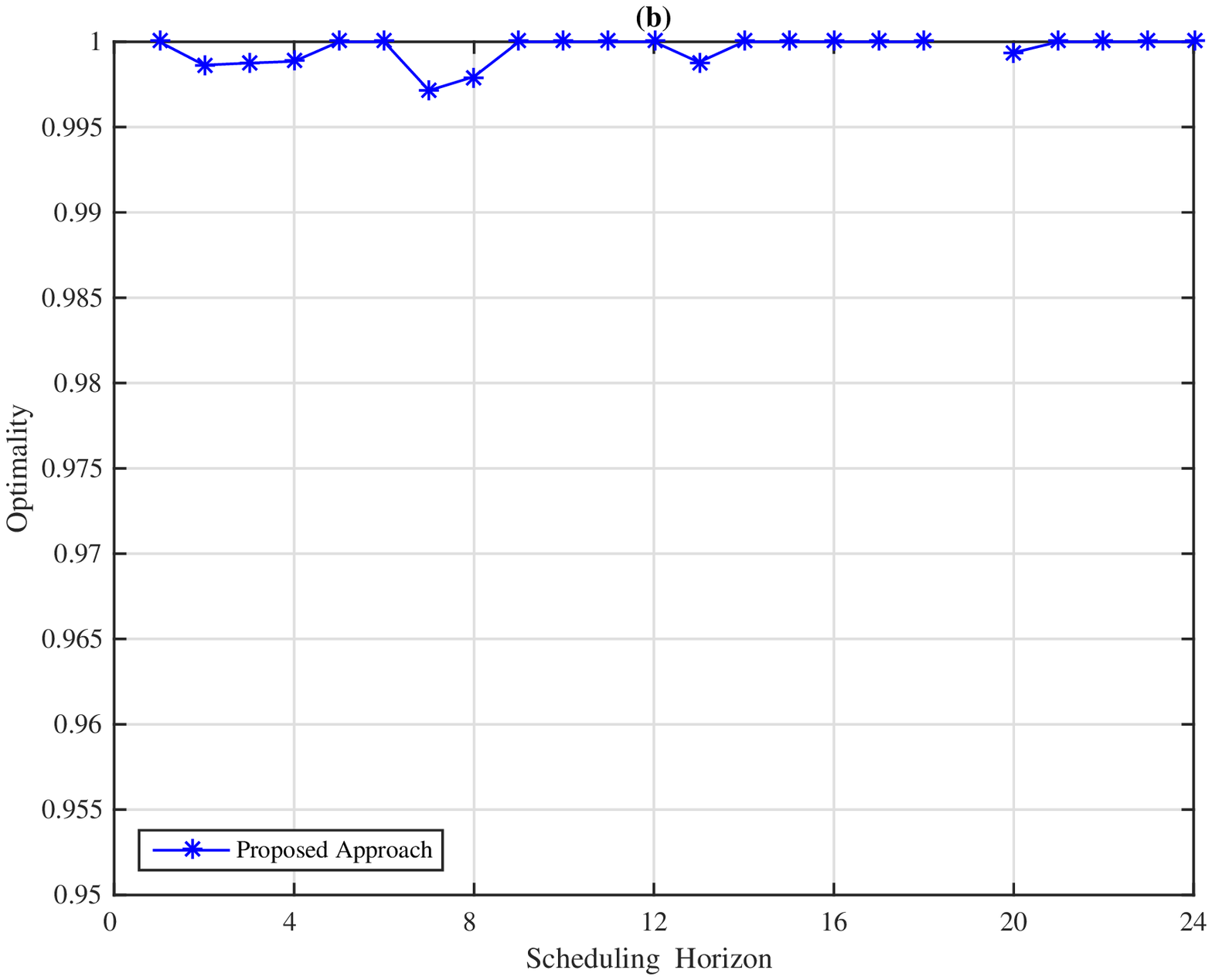}}} }
    \
      
    \vspace{-.5cm}
      \caption{The impact of increasing the optimization scheduling horizon on the performance of the proposed approach and the approach in \cite{ruiz2009pool}: (a) the computation time; (b) the optimality.}
       \vspace{-.0cm} 
 \label{fig:2}
\end{figure}

 Next, we examine the impact of changing the scheduling horizon.  To allow the competing MILP approach in \cite{ruiz2009pool} to converge in a timely manner, we assume that  $K = 2$,   and we instead increase the number of time slots $T$. The results are shown in Fig. \ref{fig:2}. We can see in   Fig. \ref{fig:2}(a) that,  the computation time of proposed approach grows linearly, as $T$ increases, while the computation time of the  MILP approach  in \cite{ruiz2009pool} grows with a significantly higher rate.  Specifically, for the case with $T=19$, the MILP approach in \cite{ruiz2009pool} does not  converge even after  running the related code for three days. In contrast, the computation time of our proposed approach is always less than 25 minutes.    Also, from Fig. \ref{fig:2}(b), our proposed approach is also always very accurate in terms of achieving the optimal profit for the strategic   producers.

\renewcommand{\thesubfigure}{}
\captionsetup{font=MyFigureFont}


\renewcommand{\thesubfigure}{}
\captionsetup{font=MyFigureFont}

\vspace{.11cm}
\subsection{The Impact of Congested Line Capacity}\label{sec:Congest}
\begin{figure}[t]
    \centering \subfigure[]{  {\scalebox{0.45} {\includegraphics{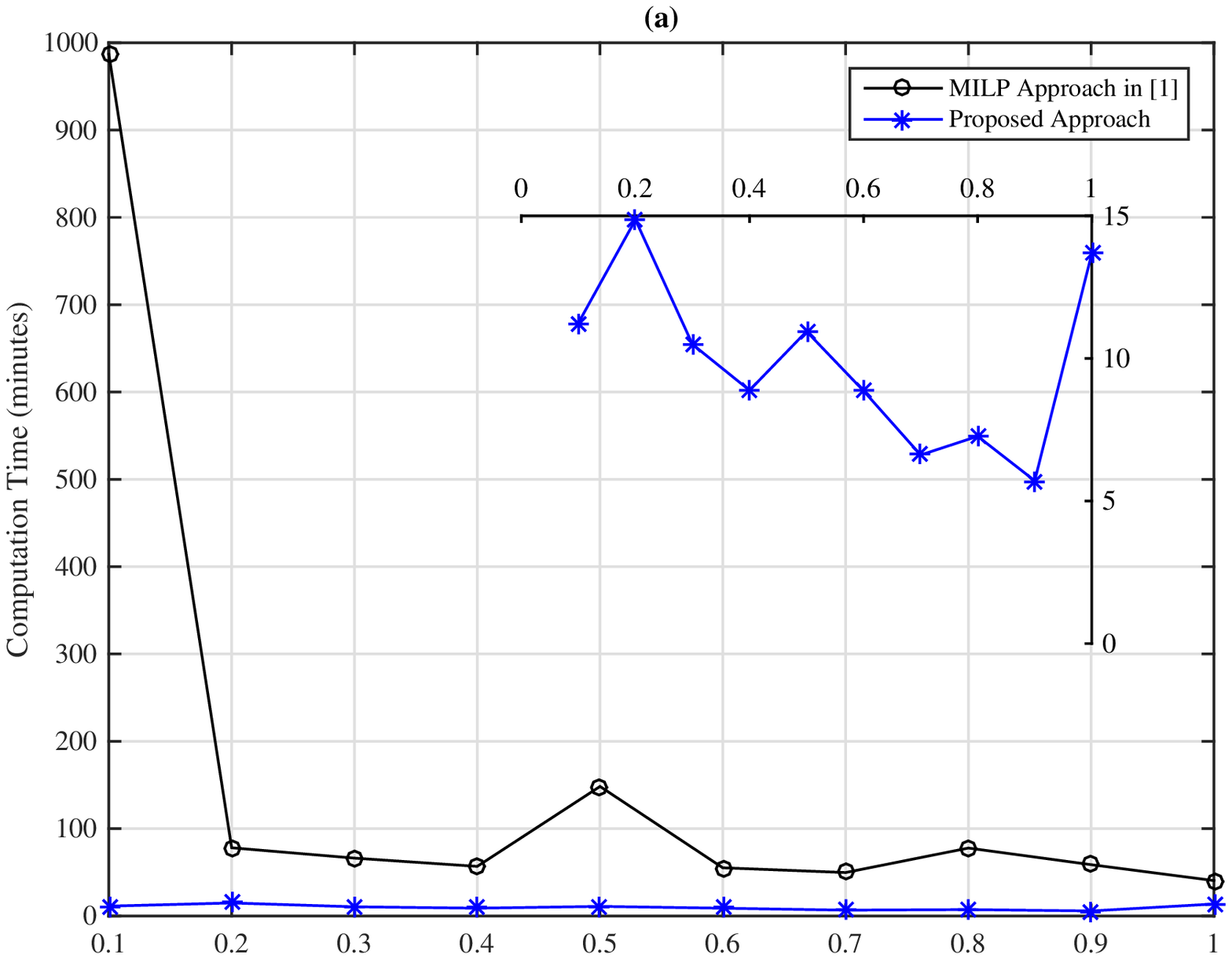}}} }
    \
      
 \vspace{-1.2cm}
 
    \centering \subfigure[]{  {\scalebox{0.45} {\includegraphics{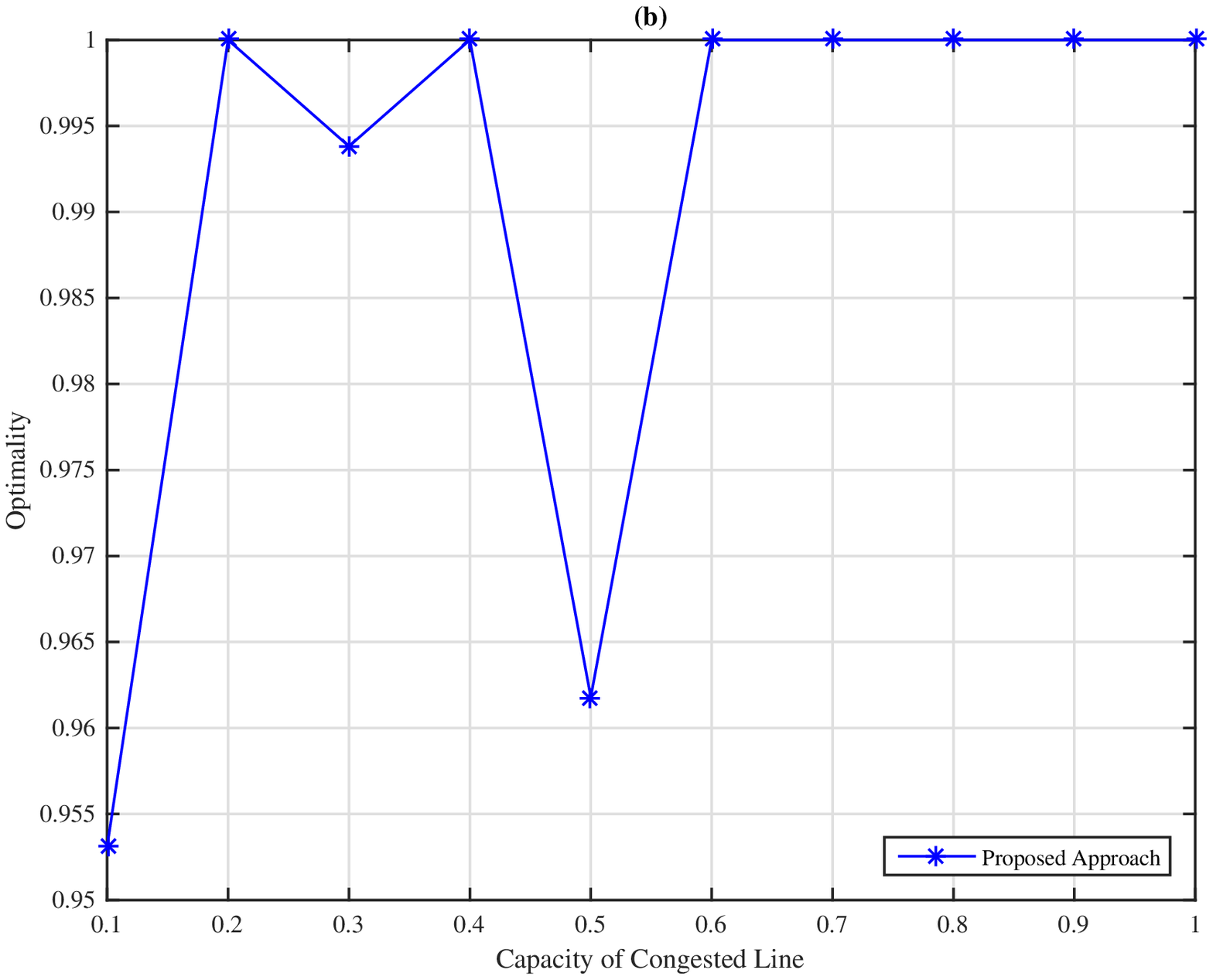}}} }
    \
      
    \vspace{-.5cm}
      \caption{The impact of changing the capacity of the congested transmission line on the performance of the proposed approach and the approach in \cite{ruiz2009pool}: (a)  the computation time; (b)  the optimality.}
  \vspace{.1cm} 
 \label{fig:3}
\end{figure}

To show that the performance of our proposed approach is not sensitive to the choice of system parameters, in this section, we examine the impact of transmission line capacity, where we set   $T=8$   and $K=3$. The results are shown in Fig. \ref{fig:3},  where we change the  capacity of transmission line 3 \cite{Mohsenian_MPEC_2015}  from 0.1 to 1.0. Again, we can see that our proposed approach is accurate and much more computationally efficient.

\renewcommand{\thesubfigure}{}
\captionsetup{font=MyFigureFont}

\renewcommand{\thesubfigure}{}
\captionsetup{font=MyFigureFont}
 \begin{figure}[h]
    \centering \subfigure[]{  {\scalebox{0.45} {\includegraphics{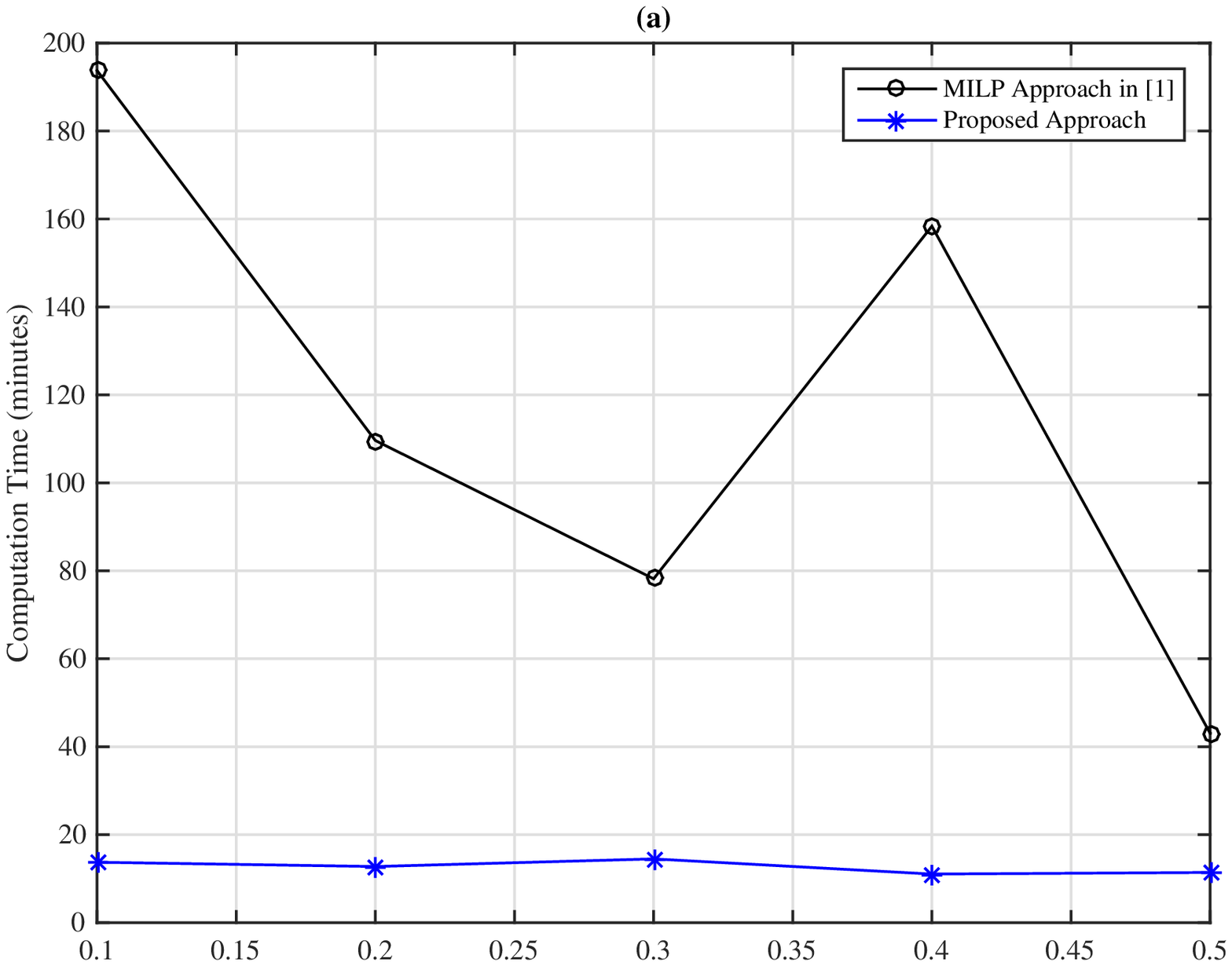}}} }
    \
      
 \vspace{-1.2cm}
 
    \centering \subfigure[]{  {\scalebox{0.45} {\includegraphics{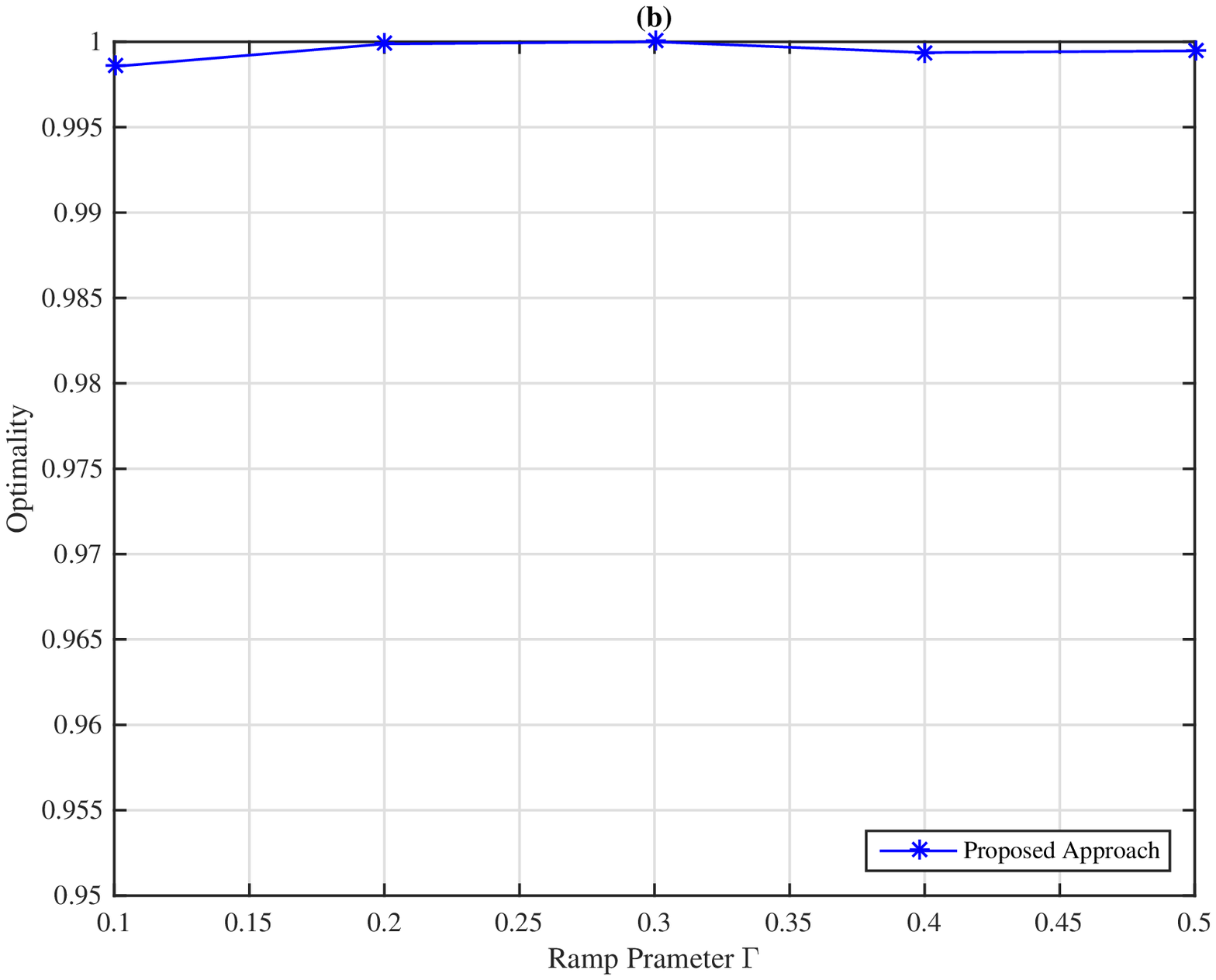}}} }
    \
      
    \vspace{-.5cm}
      \caption{ The impact of  changing  the ramp parameter $\Gamma$ on the performance of the proposed approach and the approach in \cite{ruiz2009pool}:  (a) the computation time; (b) the optimality. \color{black}} 
   \vspace{.1cm} 
 \label{fig:6}
\end{figure}

\vspace{.11cm}
\subsection{The impact of Ramp Parameter}
In this Section, the impact of the ramp parameter $\Gamma$ on the computation time as well as on  the optimality  of our proposed approach is assessed for the same  simulation setup in Section \ref{sec:Congest}, where the capacity of the  congested transmission  line is $0.2$ and the ramp parameter  $ \Gamma$ varies from $0.1$ to $0.5$. The results are shown in  Fig. \ref{fig:6}.  We can see that  our proposed approach  significantly outperforms   the MILP approach.

\vspace{.11cm}
\subsection{Comparison with other Convex Relaxation Approaches  }

In this Section,  the performance of our proposed approach is compared with that of the ones  in \cite{haghighat2014strategic} and   the SDP relaxation approaches in \cite{Fampa} and \cite{Xiaowei_2011}. The comparison is done based on the case of the IEEE 30-Bus System in Fig. \ref{fig:IEEE30bus}, where  $K=1$, and $T$ varies from 1 to 5. First and foremost, we note that    \cite{haghighat2014strategic},  \cite{Fampa} and \cite{Xiaowei_2011} do \emph{not} provide  any feasible solution to problem (\ref{formula:MPEC_KKT}). This is a common problem in many standard SDP relaxation techniques, c.f. \cite{Peng2014}. Accordingly, we can only compare the objective values under relaxation, i.e., the relaxation gap. With that in mind, we note that 
  the approach in \cite{haghighat2014strategic}   always results in an \emph{unbounded} objective value, which suggests an extremely poor performance.   The approach in \cite{Xiaowei_2011}   results in unbounded objective values for $T=1$ and $T=2$.  This approach  does not  converge for $T>2$. Therefore, the performance of the approach in \cite{Xiaowei_2011} is very poor too.  Finally, the approach in \cite{Fampa} \emph{does} converge and it \emph{is} bounded for the  cases of $T=1$ and $T=2$. This convergence is achieved after   $1239$ and $63315$ seconds, with a relaxation gap of   $5534\%$ and $3386\%$, respectively.  In contrast, once our  approach is used, the convergence times are only $22$ and $49$ seconds, and the relaxation gaps are  only $0.07\%$ and $0.15\%$, respectively. As for the cases with $T > 2$, the approach in \cite{Fampa} does \emph{not} converge.  From the above results, we can see that our proposed approach clearly outperforms  the approaches in \cite{Fampa}, \cite{haghighat2014strategic} and \cite{Xiaowei_2011}.

\color{black}

\vspace{.11cm}
\subsection{The Impact of the Number of Buses}

\begin{figure}[t]
    \centering \subfigure[]{  {\scalebox{0.45} {\includegraphics{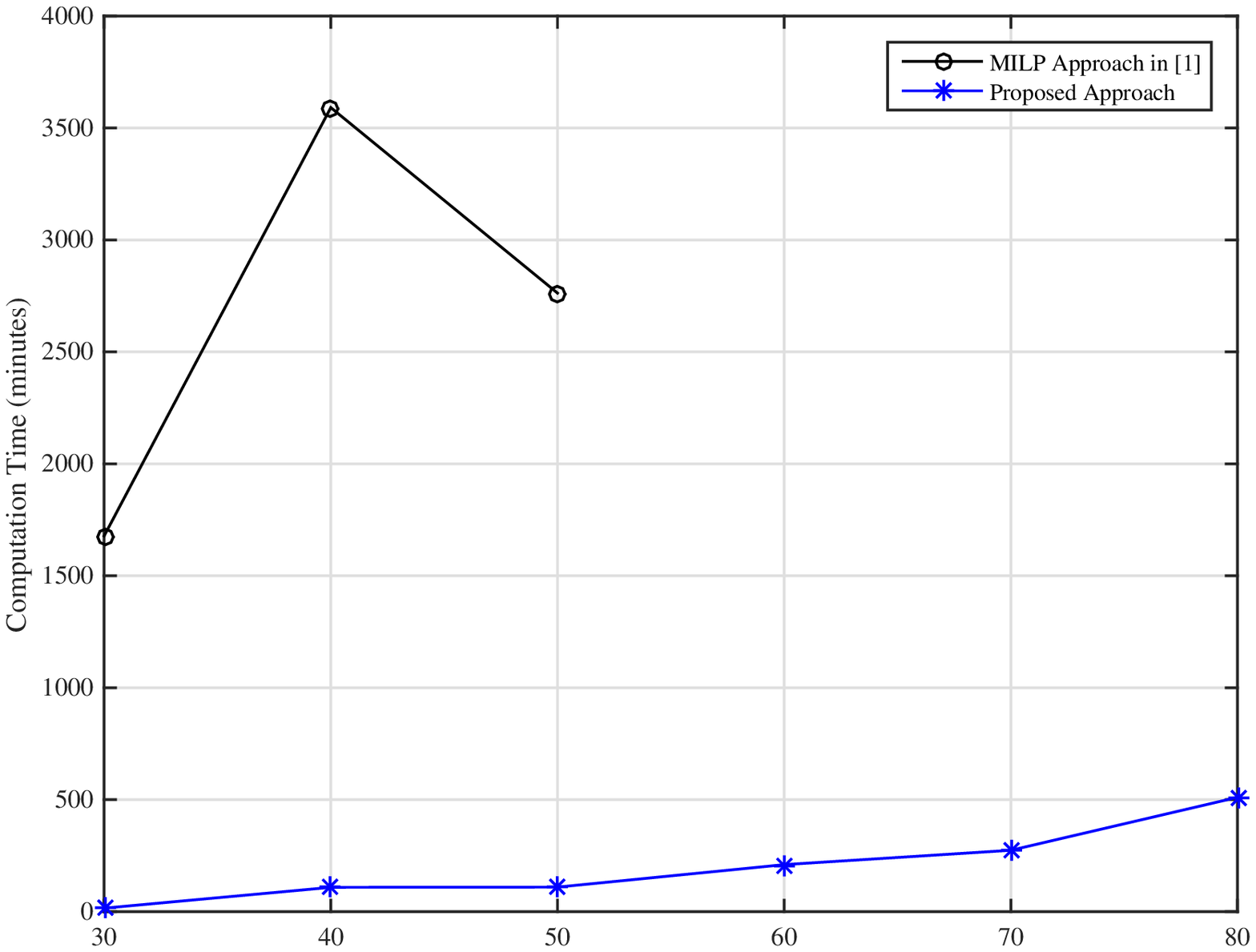}}} }
    \
      
 \vspace{-1.2cm}
 
    \centering \subfigure[]{  {\scalebox{0.45} {\includegraphics{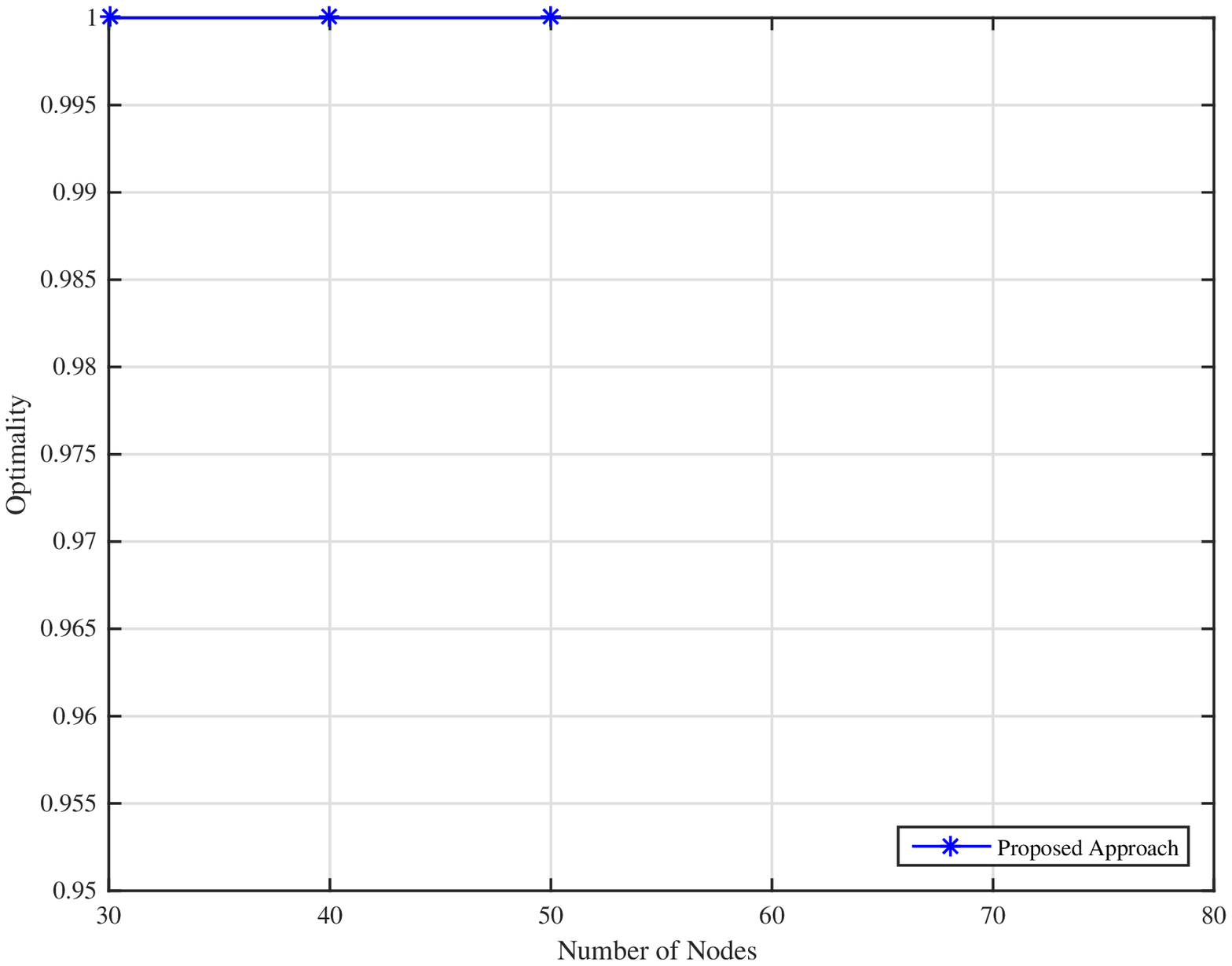}}} }
    \
      
    \vspace{-.5cm}
   \caption{ The impact of increasing the number of buses on the performance of the proposed approach and the approach in \cite{ruiz2009pool}:  (a) the computation time; (b) the optimality.  \color{black}}   
   \vspace{.1cm}
 \label{fig:15}
\end{figure}

In this Section,  the impact of the size of the  power grid on the performance of our proposed approach is assessed.   For this purpose,  several power networks are constructed by extending the number of buses, loads and generators in our base test cases  according to Table \ref{table:Bus_Number}.  
 The energy demands of the added generators are chosen such that the total added generation  is equal to the total added load. In addition,  the price bids for the added generators and  added loads are set to zero and 72 \$/MWh,   respectively.   The  line with finite capacity and the location of strategic generators are as in Section \ref{sec:SingleTimeSlot}.  Fig \ref{fig:15}(a) and Fig. \ref{fig:15}(b)  show the  computation time and  the optimality  of our proposed approach, respectively,  for the case of   $T=10$  time slots  and  $K=3$  random  scenarios.  From Fig. \ref{fig:15}(a), the computation time of our approach is much lower than  the MILP approach in \cite{ruiz2009pool}. Note that, for the power networks with 60,  70 and 80 buses, the MILP approach did not converge after three days running time.  Also, from Fig. \ref{fig:15}(b) the optimality of our approach is greater than $99\%$  for power networks with 50 buses or less. As for the cases with more than 50 buses, we simply do not know the level of optimality because we do not have a truly optimal reference for comparison. As for the networks with over 80 buses, the computation time even for our proposed approach starts growing significantly. \color{black}

\color{black}

\begin{table}[t]
\caption{Constructed Networks}
\begin{center}

\begin{tabular}{| c | c | c | c | c |  c | c | }
           
\hline \hspace{-.2cm}
  Buses & 30 & 40 & 50 & 60 & 70  & 80    \\
  \hline
  Generators  & 12 & 15 & 17 & 19 & 21   & 23  \\
    \hline
  Loads  & 16 & 21 & 26 & 31  & 35  & 41 \\
  \hline
\end{tabular}
\end{center}
\label{table:Bus_Number}
\end{table}

\vspace{.1cm}
\section{Conclusions}

A new and innovative method was proposed to solve strategic bidding problems in nodal electricity markets. Without loss of generality, we focused on the
case of strategic bidding for producers. Unlike the state-of-the-art solution approach, where the strategic bidding problem is reformulated to an MILP, the approach in this paper is based on convex programming. Therefore, in addition to its potential in achieving very accurate optimal solutions, the proposed approach is much more reliable and often computationally more tractable in solving the strategic bidding problems in power systems. For example, in a case study based on an IEEE 30-bus network with 10 random scenarios, while the state-of-the-art MILP approach does not converge even after running for about three days, our approach achieved the solution in less than twenty one minutes, running on the same computation platform. 

 While the proposed approach in this paper takes a major leap in solving strategic bidding problems in nodal electricity markets compared to the state-of-the-art MILP-based approaches, it still faces some limitations that could be addressed in future follow up studies. For example, it appears that the proposed method is well-capable of handling the increases in the number of time slots and the number of random scenarios. However, it is still not fully capable of handling the increases in the number of buses. 
Another interesting direction for future work is to obtain analytical performance bounds, i.e., on optimality and computational time, of the proposed method. \color{black}


\renewcommand{\thesubfigure}{}
\captionsetup{font=MyFigureFont}

\section*{Appendix: Proof of Theorem \ref{theorem:Optimality}}

From (\ref{formula:xstar}), the objective value of (\ref{formula:QPQC_Reduced}) at $y=y^\star$ becomes:
\begin{equation}\label{formula:Trace_F}
\begin{aligned}
&{\big(O{y^\star}+\bar{x}\big)}^TF\big(O{y^\star}+\bar{x}\big)+2f^T\big(O{y^\star}+\bar{x}\big)=\\
 & tr\bigg(
{\begin{bmatrix}  1 \\
{y^\star}
\end{bmatrix}}^T
\!{\Omega}^T
\begin{bmatrix} 0 & f^T \\ f  & F \end{bmatrix}
\Omega
\begin{bmatrix}  1 \\
{y^\star}
\end{bmatrix}
 \bigg)  \! =\! tr\left(
{\Omega}^T
\begin{bmatrix} 0 & f^T \\ f  & F \end{bmatrix} \Omega{Y^\star} \right)\! ,
\end{aligned}
\end{equation}
where the last equality is due to the fact that since $\text{Rank}(Y^\star) = 1$, $Y_{11}^\star = 1$, and (\ref{formula:xstar}) holds, we have:
\begin{equation}
Y^\star= \begin{bmatrix} 1 \\ y^\star \end{bmatrix}{ \begin{bmatrix} 1 \\ y^\star \end{bmatrix}}^T.
\end{equation}
By taking the same steps, one can show that $y=y^\star$ satisfies the constraints in problem (\ref{formula:QPQC_Reduced}).
%
Therefore, on one hand, ${y^\star}$ in (\ref{formula:xstar}) satisfies all the constraints in problem (\ref{formula:QPQC_Reduced}) and produces an objective value for problem (\ref{formula:QPQC_Reduced}) that is equal to the optimal objective value of  problem (\ref{formula:SDP_Convert2}). On the other hand, since problem (\ref{formula:SDP_Convert2}) is a convex relaxation of problem (\ref{formula:QPQC_Reduced}), its optimal objective value gives  an upper bound for the optimal objective value of  problem (\ref{formula:QPQC_Reduced}). Hence, $y^\star$ is an optimal solution for problem (\ref{formula:QPQC_Reduced}) and the relaxation gap is zero. \hfill $\blacksquare$

\vspace{0.25cm}

\bibliographystyle{IEEEtran}
\bibliography{IEEEabrv,Manuscript}

\end{document}